\documentclass[12pt, a4paper]{article}

\usepackage{graphicx}   %
\usepackage{enumerate}  %
\usepackage{latexsym}   %
\usepackage{amsmath}    %
\usepackage{amsbsy}     %
\usepackage{amsfonts}   %
\usepackage{mathrsfs}   %
\usepackage{amssymb}   %
\usepackage{array}      %
\usepackage{theorem}
\usepackage{color}
\usepackage{bm}

\usepackage{comment}

\newcommand{\jump}[1]{\ensuremath{[\![#1]\!]}}

\usepackage[dvipdfm, colorlinks=true]{hyperref}
\hypersetup{urlcolor=blue, citecolor=blue, linkcolor=blue}

 \newtheorem{definition}{Definition}
 
 \newtheorem{remark}{Remark}[section]
 \newtheorem{corollary}{Corollary}[section]
 \newtheorem{lemma}{Lemma}[section]
 
 \newtheorem{notn}{Notation}
 \numberwithin{equation}{section}

\def\bS{{\mathbb{S}}}
\def\R{{\mathbb R}}
\def\bB{{\mathbb B}}
\def\N{{\mathbb N}}

\def\L{{\mathscr L}}
\def\M{{\mathscr M}}
\def\sH{{\mathscr H}}
\def\sV{{\mathscr V}}
\def\sW{{\mathscr W}}

\def\sK{{\mathscr K}}
\def\ds{\displaystyle}
\def\ts{\textstyle}
\def\laplace{{\mit \Delta}}

\def\trace{{_{|_{\Gamma}}}}
\def\traceR{{_{|_{\R^{N-1}}}}}

\begin{document}
\thispagestyle{empty}

\pagestyle{myheadings}

\vspace*{7ex}
\begin{center}
  {\LARGE\bf {\bf
	Weak formulation for singular diffusion equation with dynamic boundary condition\footnotemark[1]
  }}
\end{center}
\vspace{3ex}
\begin{center}
{\sc Ryota Nakayashiki}\footnotemark[2]\\
e-mail: {\ttfamily nakayashiki1108@chiba-u.jp}
\\[3ex]
{\sc Ken Shirakawa}\footnotemark[3]\\
e-mail: {\ttfamily sirakawa@faculty.chiba-u.jp}
\\[7ex]
{\it Dedicated to Professor Gianni Gilardi on the occasion of his 70th birthday}
\end{center}
\vspace{5ex}

\noindent
{\bf Abstract.} In this paper, we propose a weak formulation of the singular diffusion equation subject to the dynamic boundary condition. The weak formulation is based on a reformulation method by an evolution equation including the subdifferential of a governing convex energy. Under suitable assumptions, the principal results of this study are stated in forms of Main Theorems A and B, which are respectively to verify: the adequacy of the weak formulation; the common property between the weak solutions and those in regular  problems of standard PDEs.
\\[4ex]
{\bf Key words and phrases:} singular diffusion equation, dynamic boundary condition, evolution equation, governing convex energy, Mosco-convergence, comparison principle.

\footnotetext[0]{
\noindent $\hskip-0.3cm \empty^*$\,AMS Subject Classification 
35K20, 
35K67,  
49J45. 

$\empty^\dag$\,Department of Mathematics and Informatics, Graduate School of Science, Chiba University, 1-33, Yayoi-cho, Inage-ku, Chiba, 263-8522, Japan.

$\empty^\ddag$\,Department of Mathematics, Faculty of Education, Chiba University, 1-33, Yayoi-cho, Inage-ku, Chiba, 263-8522, Japan. This author is supported by Grant-in-Aid No. 16K05224, JSPS.
}
\newpage
\section*{Introduction}
\label{sec:1}

Let $ \varepsilon > 0 $, $ 0 < T < \infty $ and $ 1 < N \in \N $ be fixed constants. Let $ \Omega \subset \R^N $ be a bounded domain with a smooth boundary $ \Gamma := \partial\Omega $, and let $ n_{\Gamma} $ be the unit outer normal to $ \Gamma $. Besides, let us denote by $ Q:= (0,T) \times \Omega $ the product space of the time interval $ (0, T) $ and the spatial domain $ \Omega $, and let us set $ \Sigma := (0, T) \times \Gamma $.

In this paper, we consider the following initial-boundary value problem of parabolic type:
\bigskip

\noindent
\begin{equation}\label{5th.eq}
\partial_{t} u - \mathrm{div}\left(\frac{Du}{|Du|} \right) = \theta(t,x) \mbox{, $ (t, x) \in Q $,}
\end{equation}
\begin{equation}\label{6th.eq}
\partial_{t} u_{\Gamma} - \varepsilon^2 {\laplace}_{\Gamma} u_{\Gamma} + (\ts{\frac{Du}{|Du|}})\trace \cdot {n}_{\Gamma} = \theta_{\Gamma}(t,x_\Gamma) \mbox{, $ (t, x_\Gamma) \in \Sigma $,}
\end{equation}
\begin{equation}\label{7th.eq}
u\trace = u_{\Gamma}(t,x_\Gamma) \mbox{, $ (t, x_\Gamma) \in \Sigma $,}
\end{equation}
\begin{equation}\label{8th.eq}
u(0,x) = u_0(x) \mbox{, $ x \in \Omega $, and } u_{\Gamma}(0,x_{\Gamma}) = u_{\Gamma,0}(x_{\Gamma}) \mbox{, $ x_{\Gamma} \in \Gamma $,}
\end{equation}
including the singular diffusion $ -\mathrm{div}(\frac{Du}{|Du|}) $ with the normal derivative $ (\ts{\frac{Du}{|Du|}})\trace \cdot {n}_{\Gamma} $. Here, $ \theta \in L^2(0,T;L^2(\Omega)) $ and $ \theta_{\Gamma} \in L^2(0,T;L^2(\Gamma)) $ are given source terms, and $ u_0 \in L^2(\Omega) $ and $ u_{\Gamma,0} \in L^2(\Gamma) $ are given initial data. $ {\laplace}_{\Gamma} $ denotes the Laplace--Beltrami operator on the surface $ \Gamma $, and ``$ \trace $'' denotes the trace of a function on $ \Omega $. In particular, the boundary conditions \eqref{6th.eq}--\eqref{7th.eq} are collectively called {\em dynamic boundary condition,} and it consists of the part of PDE \eqref{6th.eq} on the surface $ \Gamma $, and the part of transmission condition \eqref{7th.eq} between the PDEs \eqref{5th.eq}--\eqref{6th.eq}.

The representative characteristics of \{\eqref{5th.eq}--\eqref{8th.eq}\} is in the point that this problem can be regarded as a type of transmission system, containing the Dirichlet type boundary-value problem  of singular diffusion equation \{\eqref{5th.eq},\eqref{7th.eq}\}. So, referring to the previous works \cite{ABCM,Moll05}, one can remark that:
\begin{description}
\item[$\bm{(\star)}$]the expressions of the singular terms in \eqref{5th.eq}--\eqref{6th.eq} and the transmission condition \eqref{7th.eq} are practically meaningless, and for the treatments in rigorous mathematics, these must be prescribed in a weak variational sense, based on the spatial regularity in the space $ BV(\Omega) $ of functions of bounded variations. 
\end{description}

To answer the remark $\bm{(\star)}$, we here adopt an idea to put: 
\begin{equation*}
\left\{ ~ \parbox{9cm}{
$ U_0 := [u_0,u_{\Gamma,0}] $ in $ \sH := L^2(\Omega) \times L^2(\Gamma) $,
\\[1ex]
$ U := [u,u_{\Gamma}] $ and $ \Theta := [\theta,\theta_{\Gamma}] $ in $ L^2(0,T;\sH) $,
}\right.
\end{equation*}
and to reformulate the transmission system $ \{ $\eqref{5th.eq}--\eqref{8th.eq}$ \} $ to the Cauchy problem of an evolution equation:
\begin{equation}\label{CP_*}
\left\{ ~ \parbox{9cm}{
$ U'(t) + \partial\Phi_{*}(U(t)) \ni \Theta(t) $ in $ \sH $, $ t \in (0,T) $,
\\[1ex]
$ U(0) = U_0 $ in $ \sH $;
}\right.
\end{equation}
which is governed by the subdifferential $ \partial \Phi_* $ of the following convex function $ \Phi_* $ on $ \sH $:
\begin{equation}\label{phi_*}
\begin{array}{rcl}
W & = & [w,w_{\Gamma}] \in \sH \mapsto \Phi_*(W) = \Phi_*(w,w_{\Gamma})
\\[2ex]
& := & \left\{ \begin{array}{lll} \multicolumn{3}{l}{\ds\int_{\Omega} |Dw| + \int_{\Gamma} |w\trace - w_{\Gamma}| \, d\Gamma + \frac{\varepsilon^2}{2} \int_{\Gamma} |\nabla_{\Gamma} w_{\Gamma}|^2 \, d\Gamma,}
\\[2ex]
& \multicolumn{2}{l}{\mbox{ if $ w \in BV(\Omega) \cap L^2(\Omega) $, $ w_{\Gamma} \in H^1(\Gamma) $,}}
\\[2ex]
\infty, & \multicolumn{2}{l}{\mbox{ otherwise;}}
\end{array}\right.
\end{array}
\end{equation}
where $ \int_{\Omega} |Dw| $ denotes the total variation of $ w \in BV(\Omega) \cap L^2(\Omega) $, and $ \nabla_{\Gamma} $ and $ d\Gamma $ denote the surface gradient and the area element on $ \Gamma $, respectively.
Besides, we simply denote by $ \sW $ the effective domain of $ \Phi_* $, i.e. 
\begin{equation*}
\sW := (BV(\Omega) \cap L^2(\Omega)) \times H^1(\Gamma), 
\end{equation*}
and we propose to define a {\em weak solution}, i.e. the solution to a {\em weak formulation} to the system \{\eqref{5th.eq}--\eqref{8th.eq}\}, as follows. 
\begin{definition}\label{DefOfSol}
\begin{em}
	A pair of functions $ [u, u_\Gamma] \in L^2(0, T; \sH) $ is called a \textit{weak solution} to \{\eqref{5th.eq}--\eqref{8th.eq}\}, iff. $ u \in W^{1,2}(0,T; L^2(\Omega)) $, $ |Du|(\Omega) \in L^{\infty}(0,T) $, $ u_{\Gamma} \in W^{1, 2}(0, T; L^2(\Gamma)) \cap L^{\infty}(0,T;H^1(\Gamma)) $ and 
\begin{align}
&\ds \int_{\Omega} \partial_{t}u(t)(u(t) - z) \, dx + \int_{\Omega} |Du(t)| + \int_{\Gamma} |u\trace(t) - u_{\Gamma}(t)| \, d\Gamma
\nonumber\\[1ex]
 +&\int_{\Gamma} \partial_{t}u_{\Gamma}(t)(u_{\Gamma}(t)- z_{\Gamma}) \, d\Gamma + \varepsilon^2 \int_{\Gamma} \nabla_{\Gamma} u_{\Gamma}(t) \cdot \nabla_{\Gamma}(u_{\Gamma}(t) - z_{\Gamma}) \, d\Gamma
\nonumber\\[1ex]
&\ds \leq \int_{\Omega} |Dz| + \int_{\Gamma} |z\trace - z_{\Gamma}| \, d\Gamma \mbox{, \ for any $ [z,z_{\Gamma}] \in \sW $.}\label{varIneq00}
\end{align}
\end{em}
\end{definition}

As a natural consequence, the above Definition \ref{DefOfSol} will raise some issues concerned with:
\begin{description}
\item[(q1)]the adequacy of Definition \ref{DefOfSol} as the variational characterization for the singular transmission system \{\eqref{5th.eq}--\eqref{8th.eq}\};
\item[(q2)]the exemplification of fine properties which sustain common properties between our weak solutions and the solutions to  regular transmission systems via the standard dynamic boundary conditions. 
\end{description}

In the issue (q1), it will be essential to ensure that:
\begin{description}
\item[$\bm{(\star\star)}$]the Cauchy problem \eqref{CP_*} can be said as an invariant formulation to define the weak solution, i.e. the finding formulation is well-established, if various approximation approaches are applied by using many kinds of relaxation methods, with any convergent orders of the relaxation arguments. 
\end{description}
Then, it will be recommended that some of such relaxation methods are involved in the numerical approaches to our singular system.
\bigskip

In view of this, we consider the following regular transmission system via the standard dynamic boundary condition:
\noindent
\begin{equation}\label{1st.eq}
\partial_{t} u - \mathrm{div}\left(\nabla f_{\delta}(\nabla u) + \kappa^2 \nabla u \right) = \theta(t,x) \mbox{, $ (t, x) \in Q $,}
\end{equation}
\begin{equation}\label{2nd.eq}
\partial_{t} u_{\Gamma} - \varepsilon^2 {\laplace}_{\Gamma} u_{\Gamma} + (\nabla f_{\delta} (\nabla u) + \kappa^2 \nabla u)\trace \cdot {n}_{\Gamma} = \theta_{\Gamma}(t,x_\Gamma) \mbox{, $ (t, x_\Gamma) \in \Sigma $,}
\end{equation}
\begin{equation}\label{3rd.eq}
u_{|_{\Gamma}} = u_{\Gamma}(t,x_\Gamma) \mbox{, $ (t, x_\Gamma) \in \Sigma $,}
\end{equation}
\begin{equation}\label{4th.eq}
u(0,x) = u_0(x) \mbox{, $ x \in \Omega $, and } u_{\Gamma}(0,x_{\Gamma}) = u_{\Gamma,0}(x_{\Gamma}) \mbox{, $ x_{\Gamma} \in \Gamma $;}
\end{equation}
as relaxed versions of \{\eqref{5th.eq}--\eqref{8th.eq}\}. Here, $ \kappa > 0 $ and $ \delta > 0 $ are given constants, and $ \nabla f_{\delta} \in L^{\infty}(\R^N)^N $ is the differential (gradient) of a convex function $ f_{\delta} \in W^{1,\infty}(\R^N) $. Besides, the sequence $ \{ f_{\delta} \}_{\delta > 0} $ is supposed to converge to the Euclidean norm $ |{}\cdot{}| $, appropriately on $ \R^N $, as $ \delta \to 0 $. 

Now, by changing the setting of $ \{ f_\delta \}_{\delta > 0} $ in many ways, we can make various approximating problems that approach to \{\eqref{5th.eq}--\eqref{8th.eq}\} as $ \kappa, \delta \to 0 $. Also, we note that the variety of $ \{ f_\delta\}_{\delta > 0} $ can cover typical numerical regularizations for singular diffusions, such as regularization by hyperbola:
\begin{equation*}
\textstyle \omega \in \R^N \mapsto f_\delta(\omega) := \sqrt{|\omega|^2 +\delta^2}, \mbox{ for $ \delta > 0 $,}
\end{equation*}
and the Yosida-regularization of Euclidean norm $ |\cdot| $, e.t.c., even if the convergence of $ \{ f_\delta \}_{\delta>0} $ is restricted to the uniform sense. Incidentally, we can take form any convergent order of the coupling $ (\kappa,\delta) \to (0,0) $, up to the choices of sequences  $\{ \kappa_n \}_{n=1}^\infty \subset \{ \kappa \} $ and $ \{ \delta_n \}_{n=1}^\infty \subset \{ \delta\} $. Such wide flexibility will be reasonable to authorize our weak formulation, and this is the principal reason why we settle the relaxation system as stated in \eqref{1st.eq}--\eqref{4th.eq}.

In addition, referring to the previous relevant works, e.g. \cite{CC13,CGNS1X,CF1,CS,GGM08}, we can see that each approximating problem \{\eqref{1st.eq}--\eqref{4th.eq}\} is equivalent to the Cauchy problem of an evolution equation:
\begin{equation}\label{CP}
\left\{ ~ \parbox{8cm}{
$ U'(t) + \partial\Phi_{\delta}^{\kappa}(U(t)) \ni \Theta(t) $ in $ \sH $, $ t \in (0,T) $,
\\[1ex]
$ U(0) = U_0 $ in $ \sH $;
}\right.
\end{equation}
which is governed by the subdifferential $ \partial \Phi_{\delta}^{\kappa} $ of a convex function $ \Phi_{\delta}^{\kappa} : \sH \to [0,\infty] $ defined as:
\begin{equation}\label{phi_d^k}
\begin{array}{rcl}
V & = & [v,v_{\Gamma}] \in \sH \mapsto \Phi_{\delta}^{\kappa}(V) = \Phi_{\delta}^{\kappa}(v,v_{\Gamma})
\\[2ex]
& := & \left\{ \begin{array}{lll}\multicolumn{3}{l}{\ds\int_{\Omega} \left( f_{\delta}(\nabla v) + \frac{\kappa^2}{2} |\nabla v|^2 \right) \, dx + \frac{\varepsilon^2}{2} \int_{\Gamma} |\nabla_{\Gamma} v_{\Gamma}|^2 \, d\Gamma,}
\\[2ex]
& \multicolumn{2}{l}{\mbox{ if $ v \in H^1(\Omega) $, $ v_{\Gamma} \in H^1(\Gamma) $ and $ v\trace = v_{\Gamma} $ in $ L^2(\Gamma) $,}}
\\[2ex]
\infty, & \multicolumn{2}{l}{\mbox{ otherwise.}}
\end{array}\right.
\end{array}
\end{equation}

Hence, for the verification of (q1), it would be effective to observe the continuous dependence between the Cauchy problems \eqref{CP_*} and  \eqref{CP}, as $ \kappa, \delta \to 0 $, for every regularizations $ \{ f_\delta \}_{\delta > 0} $.  
Furthermore, on account of the general theories of nonlinear evolution equations and their variational convergence \cite{Attouch, Barbu,Brezis,Kenmochi81}, the essence of (q1) can be reduced as follows.
\begin{description}
\item[(A)]An issue to verify that the convex function $ \Phi_* $ on $ \sH $, given in \eqref{phi_*}, is a limit of various sequences of relaxed convex functions $ \{\Phi_{\delta}^{\kappa}\}_{\kappa, \delta > 0} $ on $ \sH $, in the sense of Mosco \cite{Mosco}, as $ \kappa, \delta \to 0 $.
\end{description}

In the meantime, for the issue (q2), we focus on the comparison principle for the weak solutions to \{\eqref{5th.eq}--\eqref{8th.eq}\}, stated as follows.
\begin{description}
	\item[(B)]If $ [u_0^k, u_{\Gamma, 0}^k] \in \sW $ and $ [\theta^k, \theta_\Gamma^k] \in L^2(0, T; \sH) $, for $ k = 1, 2 $, and
\vspace{-0.2cm}
\begin{equation*}
\left\{ ~ \parbox{7.5cm}{
$ u_0^1 \leq u_0^2 $ a.e. in $ \Omega $, $ \theta^1 \leq \theta^2 $ a.e. in $ Q $,
\\[0ex]
$ u_{\Gamma, 0}^1 \leq u_{\Gamma, 0}^2 $ a.e. on $ \Gamma $, $ \theta_\Gamma^1 \leq \theta_\Gamma^2 $ a.e. on $ \Sigma $,
} \right.
\vspace{-0.2cm}
\end{equation*}
then, it holds that:
\begin{center}
$ u^1 \leq u^2 $ a.e. in $ Q $, and $ u_{\Gamma}^1 \leq u_{\Gamma}^2  $ a.e. on $ \Sigma $, 
\end{center}
where for every $ k = 1, 2 $, $ [u^k, u_{\Gamma}^k] \in L^2(0, T; \sH) $ is a solution to \{\eqref{5th.eq}--\eqref{8th.eq}\} in the case when $ [u_0, u_{\Gamma, 0}] = [u_0^k, u_{\Gamma, 0}^k] $ and $ [\theta, \theta_\Gamma] = [\theta^k, \theta_\Gamma^k] $. 
\end{description} 
Indeed, in regular systems like \{\eqref{1st.eq}--\eqref{4th.eq}\}, the property kindred to (B) can be verified, immediately, by applying usual methods as in \cite{ABCM, Barbu, Brezis, Kenmochi81, Moll05}. But in our study, the issue of comparison principle (B) will be delicate, because the boundary integral $ \int_\Gamma |w\trace -w_\Gamma| \, d \Gamma $ as in \eqref{phi_*} will bring non-trivial interaction between the unknowns $ u $ and $ u_\Gamma $ in the transmission system  \{\eqref{5th.eq}--\eqref{8th.eq}\}. 
\medskip

In view of these, the discussions for the above (A) and (B) are developed in accordance with the following contents. In Section 1, we prepare preliminaries of this study, and in Section 2, we state the results of this paper. The principal part of our results are stated as Main Theorems A and B, and these correspond to the issues (A) and (B), respectively. Then, the continuous dependence between Cauchy problems \eqref{CP_*} and \eqref{CP} will be mentioned as a Corollary of Main Theorem A. The results are proved through the following Sections 3 and 4, which are assigned to the preparation of Key-Lemmas, and to the body of the proofs of Main Theorems and the corollary, respectively. Furthermore, in the final Section 5, we mention about an advanced issue as the future prospective of this study.

\section{Preliminaries}
\label{sec:2}

In this section, we outline some basic matters, as preliminaries of our study.
\begin{notn}[Notations in real analysis]\label{Note00}
\begin{em}
For arbitrary $ a, b \in [-\infty, \infty] $, we define:
\begin{center}
$ a \vee b := \max \{ a, b \} $ and $ a \wedge b := \min \{ a, b \} $;
\end{center}
and in particular, we write $ [a]^+ := a \vee 0 $ and $ [b]^- := -(0 \wedge b) $.
\pagebreak 

Let $ d \in \N $ be any fixed dimension. Then, we simply denote by $ |x| $ and $ x \cdot y $ the Euclidean norm of $ x \in \R^d $ and the standard scalar product of  $ x, y \in \R^d $, respectively.
Also, we denote by $ \bB^d $ and $ \bS^{d -1} $ the $ d $-dimensional unit open ball centered at the origin, and its boundary, respectively, i.e.:
\begin{equation*}
\bB^d := \left\{ \begin{array}{l|l}
x \in \R^d ~ & ~ |x| < 1
\end{array} \right\} \mbox{ and \ } \bS^{d -1} := \left\{ \begin{array}{l|l}
x \in \R^d ~ & ~ |x| = 1
\end{array} \right\}.
\end{equation*}
In particular, when $ d > 1 $, we let:
\begin{equation*}
\left\{ ~ \parbox{10cm}{
$ x \vee y := \bigl[ x_1 \vee y_1, \dots, x_d \vee y_d \bigl] $, $ x \wedge y := \bigl[x_1 \wedge y_1, \dots, x_d \wedge y_d \bigr] $,
\\[1ex]
$ [x]^+ := \bigl[ [x_1]^+, \dots, [x_d]^+ \bigr] $ and $ [y]^- := \bigl[ [y_1]^- , \dots, [y_d]^- \bigr] $,
} \right. \mbox{ \ for all $ x, y \in \R^d $.}
\end{equation*}
Besides, we often describe a  $ d $-dimensional vector $ x = [x_1, \dots, x_d] \in \R^d $ as $ x = [\tilde{x}, x_d] $ by putting $ \tilde{x} = [x_1, \dots, x_{d -1}] \in \R^{d -1} $. As well as, we describe the gradient $ \nabla = [\partial_1, \dots, \partial_d] $ as $ \nabla = [\tilde{\nabla}, \partial_d] $ by putting $ \tilde{\nabla} = [\partial_1, \dots, \partial_{d -1}] $, and additionally, we describe $ \nabla_x $, $ \partial_t$, $ \partial_{x_i} $, $ i=1,\dots,d $, and so on, when we need to specify the variables of differentials. 
\end{em}
\end{notn}
\begin{notn}[Notations of functional analysis]\label{NoteFuncSp}
\begin{em}
For an abstract Banach space $ X $, we denote by $ |{}\cdot{}|_X $ the norm of $ X $, and denote by $ {}_{X^*} \langle{}\cdot{}, {}\cdot{}\rangle_X $ the duality pairing between $ X $ and the dual space $ X^* $ of $ X $. 
In particular, when $ X $ is a Hilbert space, we denote by $ ({}\cdot{},{}\cdot{})_X $ the inner product in $ X $. 
\end{em}
\end{notn}
\begin{notn}[Notations in convex analysis]\label{NoteConvex}
\begin{em}
Let $ X $ be an abstract real Hilbert space. 

For any closed and convex set $ \mathscr{C} \subset X $, we denote by $ \pi_\mathscr{C} : X \to \mathscr{C} $ the orthogonal projection onto $ \mathscr{C} $. 

For any proper lower semi-continuous (l.s.c. from now on) and convex function $ \Psi $ defined on $ X $, we denote by $ D(\Psi) $ its effective domain, and denote by $ \partial \Psi $ its subdifferential. The subdifferential $ \partial \Psi $ is a set-valued map corresponding to a weak differential of $ \Psi $, and it has a maximal monotone graph in the product space $ X \times X $. More precisely, for each $ z_0 \in X $, the value $ \partial \Psi(z_0) $ is defined as a set of all elements $ z_0^* \in X $ which satisfy the following variational inequality:
\begin{equation*}
(z_0^*, z -z_0)_X \leq \Psi(z) -\Psi(z_0) \mbox{, for any $ z \in D(\Psi) $.}
\end{equation*}
The set \ 
$ D(\partial \Psi) := \{ z \in X \,|\, \partial \Psi(z) \ne \emptyset \} $ \ 
is called the domain of $ \partial \Psi $, and the notation ``$ [z_0, z_0^*] \in \partial \Psi $ in $ X \times X $\,'' is often rephrased as ``$ z_0^* \in \partial \Psi(z_0) $ in $ X $ with $ z_0 \in D(\partial \Psi) $'', by identifying the operator $ \partial \Psi $ with its graph in $ X \times X $.
\end{em}
\end{notn}

On this basis, we here recall the notion of Mosco-convergence for sequences of convex functions.
\begin{definition}[Mosco-convergence: cf. \cite{Mosco}]\label{Def.Mosco}
\begin{em}
Let $ X $ be an abstract Hilbert space. Let $ \Psi : X \rightarrow (-\infty, \infty] $ be a proper l.s.c. and convex function, and let $ \{ \Psi_n \}_{n = 1}^\infty $ be a sequence of proper l.s.c. and convex functions $ \Psi_n : X \rightarrow (-\infty, \infty] $, $ n \in \N $.  Then, it is said that $ \Psi_n \to \Psi $ on $ X $, in the sense of Mosco \cite{Mosco}, as $ n \to \infty $, iff. the following two conditions are fulfilled.
\begin{description}
\item[({M1}) Lower-bound condition:]$ \varliminf_{n \to \infty} \Psi_n(\check{z}_n) \geq \Psi(\check{z}) $, if \,$ \check{z} \in X $, $ \{ \check{z}_n \}_{n = 1}^\infty \subset X $, and ~ $ \check{z}_n \to \check{z} $ weakly in $ X $ as $ n \to \infty $.
\item[({M2}) Optimality condition:]for any $ \hat{z} \in D(\Psi) $, there exists a sequence $ \{ \hat{z}_n \}_{n = 1}^\infty \subset X $ such that $ \hat{z}_n \to \hat{z} $ in $ X $ and $ \Psi_n(\hat{z}_n) \to \Psi(\hat{z}) $, as $ n \to \infty $.
\end{description}
\end{em} 
\end{definition}

\begin{notn}[Notations in basic measure theory: cf. \cite{AFP,ABM}]\label{Notemeasure}
\begin{em}
For any $ d \in \N $, the $ d $-dimensional Lebesgue measure is denoted by $ \L^d $, and unless otherwise specified, the measure theoretical phrases, such as ``a.e.'', ``\,$ dt $\,'', ``\,$ dx $\,'', and so on, are  with respect to the Lebesgue measure in each corresponding dimension. Also, in the observations on a $ C^1 $-surface $ S $, the phrase ``a.e.'' is with respect to the Hausdorff measure in each corresponding Hausdorff dimension, and the area element on $ S $ is denoted by $ dS $.

Let $ d \in \N $ be any dimension, and let $ A \subset \R^d $ be any open set. We denote by $ \M(A) $ (resp. $ \M_{\rm loc}(A) $) the space of all finite Radon measures (resp. the space of all Radon measures) on $ A $. In general, the space $ \M(A) $ (resp. $ \M_{\rm loc}(A) $) is known as the dual of the Banach space $ C_0(A) $ (resp. dual of the locally convex space $ C_{\rm c}(A) $).
\end{em}
\end{notn}

\begin{notn}[Notations in BV-theory: cf. \cite{AFP, ABM, EG, G}]\label{NoteBV}
\begin{em}
Let $ d \in \N $ be a dimension of the Euclidean space $ \R^d $, and let $ A \subset \R^d$ be an open set. A function $ u \in L^1(A) $ (resp. $ u \in L_{\rm loc}^1(A) $)  is called a function of bounded variation, or a BV-function (resp. a function of locally bounded variation, or a BV$\empty_{\rm loc}$-function) on $ A $, iff. its distributional differential $ D u $ is a finite Radon measure on $ A $ (resp. a Radon measure on $ A $), namely $ D u \in \M(A) $ (resp. $ D u \in \M_{\rm loc}(A) $).
We denote by $ BV(A) $ (resp. $ BV_{\rm loc}(A) $) the space of all BV-functions (resp. all BV$\empty_{\rm loc}$-functions) on $ A $. For any $ u \in BV(A) $, the Radon measure $ D u $ is called the variation measure of $ u $, and its total variation $ |Du| $ is called the total variation measure of $ u $. Additionally, the value $|Du|(A)$, for any $u \in BV(A)$, can be calculated as follows:
\begin{equation*}
|Du|(A) = \sup \left\{ \begin{array}{l|l}
\ds \int_{A} u \ {\rm div} \,\varphi \, dy \,& \,\varphi \in C_{\rm c}^{1}(A)^d \ \ \mbox{and}\ \ |\varphi| \le 1\ \mbox{on}\ A
\end{array}
\right\}.
\end{equation*}
The space $BV(A)$ is a Banach space, endowed with the following norm:
\begin{equation*}
|u|_{BV(A)} := |u|_{L^{1}(A)} + |D u|(A),\ \ \mbox{for any}\ u\in BV(A).
\end{equation*}
Also, $ BV(A) $ is a metric space, endowed with the following distance:
$$
[u, v] \in BV(A)^2 \mapsto |u -v|_{L^1(A)} +\left| \int_A |Du| -\int_A |Dv| \right|.
$$
The topology provided by this distance is called the {\em strict topology} of $ BV(A) $ and the convergence of sequence in the strict topology is often phrased as ``strictly in $ BV(A) $''.

In particular, if $ d > 1 $, if the open set $ A $ is bounded, and if the boundary $\partial A$ is Lipschitz, then the space $BV(A)$ is continuously embedded into $L^{d/(d -1)}(A)$ and compactly embedded into $L^{q}(A)$ for any $1 \le q < d/(d-1)$ (cf. \cite[Corollary 3.49]{AFP} or \cite[Theorem 10.1.3--10.1.4]{ABM}). Besides, there exists a (unique) bounded linear operator $ \mathcal{T}_{\partial A} : BV(A) \mapsto L^1(\partial A) $, called \em trace, \em such that $ \mathcal{T}_{\partial A} \varphi = \varphi|_{\partial A} $ on $ \partial A $ for any $ \varphi \in C^1(\overline{A}) $. Hence, in this paper, we shortly denote the value of trace $ \mathcal{T}_{\partial A} u \in L^1(\partial A) $ by $ u_{|_{\partial A}} $. 
Additionally, if $1 \le r < \infty$, then the space $C^{\infty}(\overline{A})$ is dense in $BV(A) \cap L^{r}(A)$ for the {\em intermediate convergence} (cf. \cite[Definition 10.1.3. and Theorem 10.1.2]{ABM}), i.e. for any $u \in BV(A) \cap L^{r}(A)$, there exists a sequence $\{u_{n} \}_{n = 1}^\infty \subset C^{\infty}(\overline{A})$ such that $u_{n} \to u$ in $L^{r}(A)$ and $\int_{A}|\nabla u_{n}|dx \to |Du|(A)$ as $n \to \infty$.\end{em}
\end{notn}

\begin{remark}\label{Rem.trace}
\begin{em}
(cf. \cite[Theorem 3.88]{AFP})
Let $ 1 < d \in \N $, and let $ A \subset \R^d $ be a bounded open set with a Lipschitz boundary $ \partial A $. Then, it holds that:
\begin{equation}
\int_{\partial A} u_{|_{\partial A}} \, (\psi \cdot n_{\partial A}) \, d \mathcal{H}^{d -1} = \int_A u \, {\rm div} \, \psi \, dx +\int_A \psi \cdot Du, \mbox{ for any $ \psi \in C_c^1(\R^d)^d $,}\nonumber
\end{equation}
where $ n_{\partial A} $ denotes the unit outer normal on $ \partial A $. Moreover, the trace $ \mathcal{T}_{\partial A} : BV(A) \to L^1(\partial A) $ is continuous with respect to the strict topology of $ BV(A) $. Namely, the convergence of continuous dependence holds: 
\begin{equation}\label{trace100}
\mathcal{T}_{\partial A} u_n \to \mathcal{T}_{\partial A} u \mbox{ \ as $ n \to \infty $, for $ u \in BV(A) $ and $ \{ u_n \}_{n=1}^{\infty} \subset BV(A) $,}
\end{equation}
in the topology of $ L^1(\partial A) $, if $ u_n \to u $ strictly in $ BV(A) $. However, in contrast with the traces on Sobolev spaces, it must be noted that the convergence \eqref{trace100} is not guaranteed, if $ u_n \to u $ weakly-$*$ in $ BV(A) $, and  even if we adopt any weak topology for \eqref{trace100} (including the distributional one).
\end{em}
\end{remark}

\begin{notn}[Extensions of functions: cf. \cite{AFP,ABM}]\label{Ext01}
\begin{em}
Let $ d \in \N $, let $ \mu $ be a positive measure on $ \R^d $, and let $ B \subset \R^d $ be a $ \mu $-measurable Borel set. 
For any $ \mu $-measurable function $ u : B \rightarrow \R $, we denote by $ [u]^{\rm ex} $ an extension of $ u $ over $ \R^d $. More precisely, $ [u]^{\rm ex} : \R^d \rightarrow\R $ is a Lebesgue measurable function such that $ [u]^{\rm ex} $ has an expression as a $ \mu $-measurable function on $ B $, and $ [u]^{\rm ex} = u $ $ \mu $-a.e. in $ B $. In general, the extension of $ [u]^{\rm ex} : \R^d \rightarrow\R $ is not unique, for each $ u: B \to \R $. 
\end{em}
\end{notn}
\begin{remark}\label{Rem.ext}
\begin{em}
Let $ 1 < d \in \N $, and let $ A \subset \R^d $ be a bounded open set with a $ C^1 $-boundary $ \partial A $. Then, for the extensions of functions in $ BV(A) $ and $ H^{\frac{1}{2}}(\partial A) $, we can check the following facts. 
\begin{description}
\item[(Fact\,1)](cf. \cite[Proposition 3.21]{AFP}) There exists a bounded linear operator $ \mathcal{E}_A : BV(A) \rightarrow BV(\R^d) $, such that:
\begin{itemize}
\item[--] $ \mathcal{E}_A $ maps any function $ u \in BV(A) $ to an extension $ [u]^{\rm ex} \in BV(\R^d) $;
\item[--] for any $ 1 \leq q < \infty $,  $ \mathcal{E}_A ({W^{1, q}(A)}) \subset W^{1, q}(\R^d) $, and the restriction $ \mathcal{E}_A  |_{W^{1, q}(A)} : W^{1, q}(A) \rightarrow W^{1, q}(\R^d) $ forms a bounded and linear operator with respect to the (strong-)topologies of the restricted Sobolev spaces.
\end{itemize}
\item[(Fact\,2)](cf. \cite[Theorem 5.4.1 and Proposition 5.6.3]{ABM}) There exists a bounded linear operator $ \mathcal{E}_{\partial A} : H^{\frac{1}{2}}(\partial A) \rightarrow H^1(\R^d) $, which maps any function $ \varrho \in H^{\frac{1}{2}}(\partial A) $ to an extension $ [\varrho]^{\rm ex} \in H^1(\R^d) $. 
\end{description}
\end{em}
\end{remark}

Next, we prepare the notations for the spatial domain $ \Omega $ and functions and measures on this domain. 
\begin{notn}[Notations for the spatial domain]\label{NoteOmega}
\begin{em}
Throughout this paper, let $ 1 < N \in \N $, let $ \Omega \subset \R^N $ be a bounded domain with a $ C^\infty $-boundary $ \Gamma := \partial \Omega $ and the unit outer normal $ n_\Gamma \in C^\infty(\Gamma)^N $. Besides, we suppose that $ \Omega $ and $\Gamma $ fulfill the following two conditions. 
\begin{description}
\item[{\boldmath($\Omega$0)}]There exists a small constant $ r_\Gamma > 0 $, and 
the mapping
\begin{equation*}
d_\Gamma : x \in \overline{\Omega} \mapsto \displaystyle \inf_{y \in \Gamma} |x -y| \in [0, \infty), 
\end{equation*}
forms a smooth function on the neighborhoods of $ \Gamma $:
$$
\Gamma(r) := \left\{ \begin{array}{l|l} 
x \in \Omega & d_\Gamma(x) < r
\end{array} \right\}, \mbox{ for every $ r \in (0, r_\Gamma] $.}
$$
\item[{\boldmath($\Omega$1)}]There exists a small constant $ r_* \in (0, r_\Gamma] $, and for any $ x_\Gamma \in \Gamma $ and arbitrary $ \rho, r \in (0, r_*] $, the neighborhood:
\begin{equation*}
G_{x_\Gamma}(\rho, r) := \left\{ \begin{array}{l|l} ~ 
y +x_\Gamma +\tau n_\Gamma & ~ \parbox{5cm}{
$ \tau \in (-r, r) $, $ y \in \Gamma -x_\Gamma $, and \\[0.5ex] $ \bigl| y -\bigl( y \cdot n_\Gamma(x_\Gamma) \bigr) n_\Gamma(x_\Gamma) \bigr| < \rho $
}
\end{array} \right\},
\end{equation*}
is transformed to a cylinder:
\begin{equation*}
\Pi_0(\rho, r) := \left\{ \begin{array}{l|l}
\xi = [\tilde{\xi}, \xi_N] \in \R^N & \tilde{\xi} \in \rho \bB^{N -1} \mbox{ and \,} \xi_N \in (-r, r)
\end{array} \right\},
\end{equation*}
by using a uniform $ C^\infty $-diffeomorphism $ \Xi_{x_\Gamma} : G_{x_\Gamma}(r_*, r_*) \to \Pi_0(r_*, r_*) $. Additionally, for any $ x_\Gamma \in \Gamma $, there exists a function $ \gamma_{x_{\Gamma}} \in C^\infty(r_* \overline{\bB^{N -1}}) $, a congruence transform $ \Lambda_{x_\Gamma} : \R^N \to \R^N $ and a $ C^\infty $-diffeomorphism $ H_{x_\Gamma} : \Lambda_{x_\Gamma} G_{x_\Gamma}(r_*, r_*) \to \Pi_0(r_*, r_*) $ such that:
\begin{description}
\item[{\boldmath ($\omega0$)}]$ \Xi_{x_\Gamma} = H_{x_\Gamma} \circ \Lambda_{x_\Gamma} $ as a mapping from $ G_{x_{\Gamma}}(r_*,r_*) $ onto $ \Pi_0(r_*,r_*)$;
\vspace{2ex}
\item[{\boldmath ($\omega1$)}]$ \gamma_{x_\Gamma}(0) = 0 $, and $ \nabla \gamma_{x_\Gamma}(0) = 0 $ in $ \R^{N -1} $;
\vspace{2ex}
\item[{\boldmath ($\omega2$)}]for every $ \rho, r \in (0, r_*] $, 
$$
\Lambda_{x_\Gamma}G_{x_\Gamma}(\rho, r) = Y_{x_\Gamma}(\rho, r) := \left\{ \begin{array}{l|l}
y = [\tilde{y}, y_N] \in \R^N & ~ \parbox{4.5cm}{
$ [\tilde{y}, y_N -\gamma_{x_{\Gamma}}(\tilde{y})] \in \Pi_0(\rho, r) $
}
\end{array}\right\},
$$
and in particular, 
\begin{equation*}
\Lambda_{x_\Gamma}\bigl( \Gamma \cap G_{x_\Gamma}(\rho, r) \bigr) = \left\{ \begin{array}{l|l}
y = [\tilde{y},\gamma_{x_\Gamma}(\tilde{y})] \in \R^N & ~ \parbox{1.9cm}{
$ \tilde{y} \in \rho \bB^{N -1} $
}
\end{array}\right\};
\end{equation*}
\item[{\boldmath~~($\omega3$)}]for every $ \rho, r \in (0, r_*] $, 
$$
H_{x_\Gamma} : y = [\tilde{y}, y_N] \in Y_{x_\Gamma}(\rho, r) \mapsto \xi = H_{x_\Gamma}y := [\tilde{y}, y_N -\gamma_{x_\Gamma}(\tilde{y})] \in \Pi_0(\rho, r).
$$
\end{description}
\end{description}
\end{em}
\end{notn}

\begin{remark}\label{Rem.Omega0}
\begin{em}
From ($\Omega$0), we may further suppose the following condition.
\begin{description}
\item[{\boldmath($\Omega$2)}]For any $ \sigma > 0 $, there exists a constant $ \rho_*^\sigma \in (0, r_*] $ such that:
\begin{equation*}
\begin{array}{c}
\ds \rho_*^\sigma \leq \sigma, ~ |\gamma_{x_\Gamma}|_{C^1(\rho \overline{\bB^{N -1}})} \leq \sigma \mbox{ and } 
\\[2ex]
\left\{\begin{array}{l|l}{} \Xi_{x_\Gamma}^{-1} [\tilde{\xi},\gamma_{x_\Gamma}(\tilde{\xi}) + r_*] \,&\, \tilde{\xi} \in \rho \overline{\bB^{N-1}} \end{array} \right\} \cap \overline{\Gamma(r_*/2)} = \emptyset, 
\\[2ex]
\mbox{for any $ x_\Gamma \in \Gamma $ and any $ \rho \in (0, \rho_*^\sigma] $.}
\end{array}
\end{equation*}
\end{description}
\end{em}
\end{remark}

\begin{notn}[Notations of surface-differentials]\label{bdryOp}
\begin{em}
	Under the assumption ($\Omega$0) in Notation \ref{NoteOmega}, we can define the Laplacian $ \laplace_\Gamma $ on the surface $ \Gamma $, i.e. the so-called Laplace--Beltrami operator, as the composition $ \laplace_\Gamma := {\rm div}_\Gamma \circ \nabla_\Gamma : C^\infty(\Gamma) \to C^\infty(\Gamma) $ of the {\em surface gradient:}
\begin{equation*}
	\nabla_\Gamma \varphi := \nabla [\varphi]^{\rm ex} -(\nabla {d}_\Gamma \otimes \nabla {d}_\Gamma) \nabla [\varphi]^{\rm ex}, \mbox{ for any $ \varphi \in C^\infty(\Gamma) $,}
\end{equation*}
and the {\em surface-divergence:}
\begin{equation*}
	{\rm div}_\Gamma \omega := {\rm div} [\omega]^{\rm ex} - \nabla ([\omega]^{\rm ex} \cdot \nabla d_{\Gamma}) \cdot \nabla d_{\Gamma}, \mbox{ for any $ \omega = [\omega_1, \dots, \omega_N] \in C^\infty(\Gamma)^N $.}
\end{equation*}
As is well-known (cf. \cite{SV97}), the values $ \nabla_\Gamma \varphi $ and $ {\rm div}_\Gamma \omega $ are determined independently with respect to the choices of the extensions $ [\varphi]^{\rm ex} \in C^\infty(\R^N) $ and $ [\omega]^{\rm ex} = \bigl[ [\omega_1]^{\rm ex}, \dots, [\omega_N]^{\rm ex} \bigr] \in C^{\infty}(\R)^N $, and moreover, the operator $ -\laplace_\Gamma $ can be extended to a duality map between $ H^1(\Gamma) $ and $ H^{-1}(\Gamma) $, via the following variational identity:
\begin{equation*}
\ds {}_{{}^{H^{-1}(\Gamma)}}\langle -{\mit \Delta}_\Gamma \varphi, \psi \rangle_{{}^{H^{1}(\Gamma)}} = (\nabla_\Gamma \varphi, \nabla_\Gamma \psi )_{L^2(\Gamma)^N}, \mbox{ for all $ [\varphi, \psi] \in H^1(\Gamma)^2 $.}
\end{equation*}
\end{em}
\end{notn}

Finally, we prescribe some specific notations. 
\begin{notn}\label{TV_Dirichlet}
\begin{em}
Let $ R_\Omega > 0 $ be a sufficiently large constant, such that $ \bB_\Omega := R_\Omega \bB^{N} \supset \overline{\Omega} $.
Besides, for any $ u \in BV(\Omega) $ and any $ g \in H^{\frac{1}{2}}(\Gamma) $, we denote by $ [u]_g^{\rm ex} \in BV(\bB_\Omega) \cap H^1(\bB_\Omega \setminus \overline{\Omega}) $ an extension of $ u $, provided as:
\begin{equation}\label{ext*}
x \in \R^N \mapsto [u]_g^{\rm ex}(x) := \left\{ \begin{array}{l}
u(x), \mbox{ if $ x \in \Omega $,}
\\[1ex]
[g]^{\rm ex}(x), \mbox{ if $ x \in \bB_\Omega \setminus \overline{\Omega} $,}
\end{array} \right. 
\end{equation}
with the use of an extension $ [g]^{\rm ex} \in H^1(\R^N) $ of $ g $.
\end{em}
\end{notn}
\begin{remark}\label{Rem.Du_Dirichlet}
\begin{em}
	As consequences of BV-theory (cf. \cite[Corollary 3.89]{AFP}, \cite[Example 10.2.1]{ABM} and \cite[Theorem 5.8]{EG}) and Remark \ref{Rem.ext}, we can verify the following facts.
\begin{description}
\item[{(Fact\,3)}]For any $ u \in BV(\Omega) $ and any $ g \in H^{\frac{1}{2}}(\Gamma) $, it holds that:
\begin{equation*}
\begin{array}{c}
\displaystyle 
|D[u]_g|(B) = \int_{B \cap \Omega} |Du| +\int_{B \cap \Gamma} |u\trace -g| \, d \Gamma +\int_{B \setminus \overline{\Omega}} |\nabla [g]^{\rm ex}| \, dx,
\\[2ex]
\mbox{for any Borel set $ B \subset \bB_\Omega $, and any extension $ [g]^{\rm ex} \in H^1(\R^N) $ of $ g $.}
\end{array}
\end{equation*}
\item[{(Fact\,4)}]For any $ g \in H^{\frac{1}{2}}(\Gamma)$, the functional:
\begin{equation*}
\begin{array}{rcl}
u \in L^1(\Omega) & \mapsto &\left| D[u]_g^{\rm ex} \right|(\overline{\Omega}) 
\\[2ex]
&:= &\left\{ \begin{array}{ll}
		\multicolumn{2}{l}{\ds \int_\Omega |Du| +\int_\Gamma |u\trace -g| \, d\Gamma = |D[u]_g^{\rm ex}|(\bB_\Omega) -|D[g]^{\rm ex}|(\bB_\Omega \setminus \overline{\Omega}),}
\\[2ex]
& \mbox{ \ if $ u \in BV(\Omega) $,}
\\[2ex]
\infty, & \mbox{ \ otherwise,}
\end{array} \right.
\end{array}
\end{equation*}
forms a single-valued proper l.s.c. and convex function on $ L^1(\Omega) $.
\item[{(Fact\,5) (cf. \cite{ABCM, Anzellotti, Temam})}]$|D[u_n]_g^{\rm ex}|(\overline{\Omega}) \to |D[u]_g^{\rm ex}|(\overline{\Omega}) $ as $ n \to \infty $, whenever $ \{ u_n \}_{n = 1}^\infty \\ \subset BV(\Omega) \cap L^2(\Omega) $, $ u \in BV(\Omega) \cap L^2(\Omega) $ and $ u_n \to u $ in $ L^2(\Omega) $ and strictly in $ BV(\Omega) $ as $ n \to \infty $.
\end{description}
\end{em}
\end{remark}
\begin{remark}\label{Rem.l.s.c.}
\begin{em}
From the definition \eqref{phi_*}, we easily see that $ \Phi_*$ is proper and convex. Also, the above Remark \ref{Rem.Du_Dirichlet} (Fact\,4)--(Fact\,5) lead to the lower semi-continuity of this $ \Phi_* $. In fact, taking arbitrary $ W=[w,w_\Gamma] \in \sH $ and $ \{ W_n=[w_n,w_{\Gamma,n}] \}_{n=1}^{\infty} \subset\sW $, such that:
\begin{equation*}
W_n=[w_n,w_{\Gamma,n}] \to W=[w,w_\Gamma] \mbox{ \ in $ \sH $, as $ n \to \infty $,}
\end{equation*}
we immediately see from Remark \ref{Rem.Du_Dirichlet} (Fact\,5) that:
\begin{align*}
\varliminf_{n\to\infty} \Phi_*(W_n) &\geq \varliminf_{n \to \infty} |D[w_n]_{w_\Gamma}^{\rm ex}|(\overline{\Omega}) - \lim_{n \to \infty} \int_\Gamma |w_\Gamma\trace - w_\Gamma|\,d\Gamma
\\
& \hspace{1.5cm} +\frac{\varepsilon^2}{2} \varliminf_{n\to\infty} \int_\Gamma |\nabla_\Gamma w_{\Gamma,n}|^2\,d\Gamma
\\
& \geq |D[u]_{w_\Gamma}^{\rm ex}|(\overline{\Omega}) + \frac{\varepsilon^2}{2} \int_{\Gamma} |\nabla_\Gamma w_\Gamma|^2\,d\Gamma = \Phi_*(W).
\end{align*}
\end{em}
\end{remark}

\section{The results of this paper}

First, we prescribe, anew, the product Hilbert space 
$ \mathscr{H} := L^2(\Omega) \times L^2(\Gamma) $, with the inner product:
\begin{equation*}
\begin{array}{c}
\left( [z^1, z_\Gamma^1], [z^2, z_\Gamma^2] \right)_{\mathscr{H}} := (z^1, z^2)_{L^2(\Omega)} +(z_\Gamma^1, z_\Gamma^2)_{L^2(\Gamma)}, 
\\[2ex]
\mbox{for all $ [z^k, z_\Gamma^k] $, $ k = 1, 2 $.}
\end{array}
\end{equation*}
As is mentioned in Introduction, the Hilbert space $ \mathscr{H} $ is to be the base-space of the convex functions as in \eqref{phi_*} and \eqref{phi_d^k}, and the Cauchy problems \eqref{CP_*} and \eqref{CP}. 
Also, let $ \sW := (BV(\Omega) \cap L^2(\Omega)) \times H^1(\Gamma) $ be the effective domain of the convex function $ \Phi_* $, given in \eqref{phi_*}, and let $ \sV $ be a closed linear subspace in the product Hilbert space $ H^1(\Omega) \times H^1(\Gamma) $, defined as:
\begin{equation*}
\sV := \left\{\begin{array}{l|l} [v,v_{\Gamma}] \in \sH & ~ \parbox{4.25cm}{$ v \in H^1(\Omega) $, $ v_{\Gamma} \in H^1(\Gamma) $ \\ and $ v\trace = v_{\Gamma} $ a.e. on $ \Gamma $} \end{array}\right\}.
\end{equation*}

Next, we prescribe the assumptions in our study.
\begin{description}
\item[(A0)]$ \varepsilon > 0 $ is a fixed constant, and $ \delta > 0 $ and $ \kappa > 0 $ are given constants. Besides, $ 1 <  N \in \N $ is a fixed constant, and $ \Omega \subset \R^N $ is a bounded domain with a smooth boundary $ \Gamma := \partial \Omega $ and the unit outer normal $ n_\Gamma $, that fulfills the conditions ($\Omega$0)--($\Omega$1) in Notation \ref{NoteOmega}.
\item[(A1)]$ \{ f_{\delta} \}_{0<\delta \leq 1} \subset W^{1,\infty}(\R^N) $ is a sequence of convex functions such that
\begin{equation*}
\begin{array}{c}
f_{\delta}(0) = 0 \mbox{ and } f_{\delta}(\omega) \geq 0 \mbox{, for any $ 0< \delta \leq 1 $ and any $ \omega \in \R^N $,}
\\[1ex]
\mbox{ and } f_{\delta} \to |{}\cdot{}|(=|{}\cdot{}|_{\R^N}) \mbox{, uniformly on $ \R^N $, as $ \delta \to 0 $.}
\end{array}
\end{equation*} 
\end{description}

\begin{remark}\label{Rem.Assume}
\begin{em}
The assumptions (A0)--(A1) cover the setting of $ \{ f_{\delta}\}_{\delta>0} = \{ |{}\cdot{}| \} $, and this setting is just the case that was mainly dealt with in the previous work \cite{CGNS1X}.
\end{em}
\end{remark}

Now, the results of this paper are stated as follows.
\paragraph{Main Theorem A (Mosco-convergence).}{\em 
Under (A1)--(A0), let $ \Phi_* : \mathscr{H} \to [0, \infty] $ be the functional given in \eqref{phi_*}, and for every $ \delta > 0 $ and $ \kappa > 0 $, let $ \Phi_\delta^\kappa : \mathscr{H} \to [0, \infty] $ be the proper l.s.c. and convex function given in \eqref{phi_d^k}. Then, 
for every sequences $ \{ \delta_n \}_{n = 1}^{\infty} \subset (0,1] $ and $ \{ \kappa_n \}_{n = 1}^{\infty} \subset (0,1] $, such that:
\begin{equation}\label{MT01}
\delta_n \to 0 \mbox{ and } \kappa_n \to 0, \mbox{ as $ n \to \infty $,}
\end{equation}
the sequence $ \{ \Phi_n \}_{n = 1}^{\infty} $ of convex functions $ \Phi_n := \Phi_{\delta_n}^{\kappa_n} : \sH \to [0, \infty] $, $ n \in \N $, converges to the convex function $ \Phi_* $ on $ \sH $, in the sense of Mosco, as $ n \to \infty $. 
}
\begin{corollary}[Continuous dependence of Cauchy problems]\label{Cor.03}
Let $ 0<T<\infty $, and for every $ U_0 =[u_0,u_{\Gamma,0}] \in \sW $ and $ \Theta =[\theta,\theta_{\Gamma}] \in L^2(0,T;\sH) $, let $ U = [u,u_\Gamma] \in L^2(0,T;\sH) $ be the solution to \eqref{CP_*}. Also, for every $ n \in \N $, $ U_0^n := [u_0^n, u_{\Gamma,0}^n] \in \sV $, and $ \Theta^n := [\theta^n,\theta_{\Gamma}^n] \in L^2(0,T;\sH) $, let $ U^n := [u^n,u_{\Gamma}^n] \in W^{1,2}(0,T;\sH) \cap L^{\infty}(0,T;\sV) $ be the solution to \eqref{CP} in the case when $ \delta = \delta_n $ and $ \kappa = \kappa_n $, i.e.:
\begin{equation*}
\left\{ ~ \parbox{10cm}{
$ (U^n)'(t) + \partial \Phi_n(U^n(t)) \ni \Theta^n(t) $ in $ \sH $, a.e. $ t \in (0,T) $,
\\[1ex]
$ U^n(0) = U_0^n $ in $ \sH $.
}\right.
\end{equation*}
On this basis, let us assume that:
\begin{equation*}
U_0^n \to U_0 \mbox{ in $ \sH $ and } \Theta^n \to \Theta \mbox{ in $ L^2(0,T;\sH) $, with \eqref{MT01}.}
\end{equation*}
Then, the sequence $ \{ U^n = [u^n,u_{\Gamma}^n] \}_{n=1}^{\infty} $ converges to $ U = [u,u_{\Gamma}] $ in the sense that:
\begin{equation*}
\begin{array}{ll}
U^n \to U & \mbox{ in $ C([0,T];\sH) $,}
\\[0.5ex]
               &\mbox{ weakly in $ W^{1,2}(0,T;\sH) $, as $ n \to \infty $,}
\end{array}
\end{equation*}
and
\begin{equation*}
\int_0^T \Phi_n(U^n(t)) \, dt \to \int_0^T \Phi_*(U(t)) \, dt, \mbox{ as $ n \to \infty $.}
\end{equation*}
\end{corollary}

\paragraph{Main Theorem B (Comparison principle).}{\em For every $ k=1,2 $, let $ [u_0^k,u_{\Gamma,0}^k] \in \sW $ be given initial data, let $ [\theta^k,\theta_\Gamma^k] \in L^2(0,T;\sH) $ be a given source term, and let $ [u^k,u_\Gamma^k] \in L^2(0, T; \sH) $ be a weak solution to \{\eqref{1st.eq}--\eqref{4th.eq}\} in the case when $ [u_0,u_{\Gamma,0}] = [u_0^k,u_{\Gamma,0}^k] $ and $ [\theta,\theta_\Gamma]=[\theta^k,\theta_\Gamma^k] $. Then, it holds that:
\begin{align}
&\bigl| [u^1-u^2]^+(t)\bigr|_{L^2(\Omega)}^2 + \bigl|[u_\Gamma^1 - u_\Gamma^2]^+(t)\bigr|_{L^2(\Gamma)}^2\nonumber
\\
\leq & {\mathrm e}^t \,\bigl( 
\bigl|[u_0^1-u_0^2]^+\bigr|_{L^2(\Omega)}^2 + \bigl|[u_{\Gamma,0}^1-u_{\Gamma,0}^2]^+\bigr|_{L^2(\Gamma)}^2 \bigr)
\nonumber 
\\
& +\int_0^t e^{t -\tau} \bigl( \bigl|[\theta^1 -\theta^2]^+(\tau) \bigr|_{L^2(\Omega)}^2 + \bigl|[\theta_\Gamma^1 -\theta_\Gamma^2]^+(\tau) \bigr|_{L^2(\Gamma)}^2 
\bigr) \, d\tau, \label{ineqB}
\\
& \hspace{3cm}\mbox{for all $ t \in [0, T] $.}
\nonumber
\end{align}
}
\begin{remark}\label{remMainB}
\begin{em}
	In Main Theorem B, we can suppose the well-posedness for the weak formulation \eqref{varIneq00}, because the Definition \ref{DefOfSol} lets the well-posedness be just a straightforward consequence of the general theory of nonlinear evolution equations \cite{Barbu,Brezis,Kenmochi81}. Also, we note that the comparison principle (B), mentioned in Introduction, is immediately deduced from the inequality \eqref{ineqB}. 
\end{em}
\end{remark}
\section{Key-Lemmas}

In Main Theorem A, the keypoint is in the construction method of approximating sequences for BV-functions, which is stated in the following Key-Lemma A. 
\paragraph{Key-Lemma A.}{
{\em
For any $ [\hat{u}, \hat{u}_{\Gamma}] \in \sW $, there exists a sequence $ \{ \hat{u}_\ell \}_{\ell=1}^{\infty} \subset H^1(\Omega) $, such that:
\begin{equation}\label{Lem100-1}
\hat{u}_\ell\trace = \hat{u}_{\Gamma} \mbox{ \,in $ H^{\frac{1}{2}}(\Gamma) $, \, for any $ \ell \in \N $,}
\end{equation}
\begin{equation}\label{Lem100-2}
\begin{array}{c}
\ds\hat{u}_\ell \to \hat{u} \mbox{ in $ L^2(\Omega) $ and }
\int_{\Omega} |\nabla \hat{u}_\ell| \, dx \to \int_{\Omega} | D\hat{u} | + \int_{\Gamma} |\hat{u}\trace - \hat{u}_{\Gamma}| \, d\Gamma, 
\\
\mbox{ as $ \ell \to \infty $.}
\end{array}
\end{equation}
}}
Meanwhile, the keypoint of Main Theorem B is in the so-called {\em $ T $-monotonicity} of the subdifferential $ \partial \Phi_* $, which is stated in the following Key-Lemma B.
\paragraph{Key-Lemma B.}{
{\em
Let $ \Phi_* $ be the convex function given in \eqref{phi_*}. Then, the subdifferential $ \partial \Phi_* $ fulfills the following inequality of $ T $-monotonicity:
\begin{equation}\label{Lem40-0}
\begin{array}{c}
\bigl( U^{*, 1} -U^{*, 2}, [U^1 -U^2]^+ \bigr)_{\sH} = (u^{*, 1} -u^{*, 2}, [u^1 -u^2]^+)_{L^2(\Omega)} 
\\[1ex]
+(u_\Gamma^{*, 1} -u_{\Gamma}^{*, 2}, [u_\Gamma^1 -u_\Gamma^2]^+)_{L^2(\Gamma)} \geq 0,
\\[1.5ex]
\mbox{ for all } [U^k, U^{*, k}] = \bigl[[u^k, u_\Gamma^k], [u^{*, k}, u_\Gamma^{*, k}] \bigr] \in \partial \Phi_* \mbox{ in $ \sH \times \sH $, $ k = 1, 2 $.}
\end{array}
\end{equation}
}}

Now, before the proofs of these Key-Lemmas, we prepare some auxiliary lemmas and remarks. 
\begin{lemma}\label{LemExt01} Let $ \R_+^N $ be the upper half-space of $ \R^N $, i.e.: 
\begin{equation*}
\R_+^N := \left\{ \begin{array}{l|l} [\tilde{\xi},\xi_N] \in \R^N & \tilde{\xi} \in \R^{N-1} \mbox{ and \ } \xi_N >0 \end{array} \right\}. 
\end{equation*}
Then, for any $ \varpi \in H^1(\R^{N -1}) \cap BV(\R^{N -1}) $, there exists a sequence $ \{ \jump{\varpi}_r^{\rm ex} \}_{r > 0} \subset H^1(\R_+^N) \cap BV(\R_+^N) $, and for any $ \tau > 0 $, there exists a small constant $ r_{\varpi}^\tau \in (0, r_*] $, such that: 
\begin{equation}\label{Lem10-0}
\begin{array}{c}
\ds r_{\varpi}^\tau \leq \tau \mbox{ \ and \ } \jump{\varpi}_r^{\rm ex}(\tilde{\xi}, \xi_N) = 0 \mbox{, \ for any $ r \in (0, r_{\varpi}^\tau] $}
\\[1ex]
\mbox{ and a.e. $ [\tilde{\xi}, \xi_N] \in \R_+^{N} $, satisfying $ \xi_N > r $,}
\end{array}
\end{equation}
\begin{equation}\label{Lem10-2}
\jump{\varpi}_r^{\rm ex}\traceR = \varpi \mbox{ in $ H^{\frac{1}{2}}(\R^{N-1}) $, \ for any $ r \in (0, r_{\varpi}^\tau] $,}
\end{equation}
and
\begin{equation}\label{Lem10-1}
\begin{array}{c}
\ds \left|\jump{\varpi}_r^{\mathrm{ex}}\right|_{L^2(\R_+^N)} \leq \tau
\mbox{ and }
\ds \left|D \jump{\varpi}_r^{\mathrm{ex}} \right|(\R_+^N) \leq |\varpi|_{L^1(\R^{N-1})} + \tau, 
\\[1ex]
\mbox{ \ for any $ r \in (0, r_{\varpi}^\tau] $.}
\end{array}
\end{equation}
\end{lemma}

\paragraph{Proof of Lemma \ref{LemExt01}.}{
Let us define:
\begin{equation}\label{Lem14}
\begin{array}{c}
	\ds \jump{\varpi}_r^{\rm ex}(\xi) = \jump{\varpi}_r^{\rm ex}(\tilde{\xi},\xi_N):= [1-r^{-1}\xi_N]^+\varpi(\tilde{\xi}), 
\\[1ex]
\mbox{for a.e. $ \tilde{\xi} \in \R^{N -1} $, a.e. $ \xi_N \geq 0 $ and any $ r > 0 $.}
\end{array}
\end{equation}
Then, from the assumption $ \varpi \in H^1(\R^{N -1}) \cap BV(\R^{N -1}) $, we immediately check that $ \{ \jump{\varpi}_r^{\rm ex} \}_{r > 0} \subset H^1(\R_+^N) \cap BV(\R_+^N) $. 

On this basis, for any $ \tau > 0 $, let us take a small constant $ r_{\varpi}^\tau \in (0, r_*] $, such that:
\begin{equation}\label{Lem13}
r_{\varpi}^\tau \in (0, \tau], ~ 
\sqrt{\frac{r_{\varpi}^\tau}{3}}|\varpi|_{L^2(\R^{N-1})} < \tau \mbox{ ~ and ~ } \frac{r_{\varpi}^\tau}{2} \int_{\R^{N-1}}|\nabla \varpi| \, d\tilde{\xi} < \tau.
\end{equation}
Then, we can see the conditions \eqref{Lem10-0}--\eqref{Lem10-2} by means of \eqref{Lem14}--\eqref{Lem13} and a standard argument of the trace.
Additionally, with \eqref{Lem14}--\eqref{Lem13} in mind, we can verify the remaining \eqref{Lem10-1} as follows.
\begin{equation*}
\begin{array}{lll}
|\jump{\varpi}_{r}^{\rm ex} |_{L^2(\R_+^N)}^2 & = &\ds \int_{\R_+^N}\left|[1 - r^{-1}\xi_N]^+ \varpi(\tilde{\xi})\right|^2 \, d\xi
\\[2ex]
& = & \ds\left(\int_0^{r} (1-r^{-1}\xi_N)^2 \, d\xi_N \right)\left(\int_{\R^{N-1}} |\varpi(\tilde{\xi})|^2 \, d\tilde{\xi}\right)
\\[2ex]
& = & \ds\frac{r}{3} \, |\varpi|_{L^2(\R^{N-1})}^2 \leq \tau^2, \mbox{ for any $ r \in (0, r_{\varpi}^\tau] $,}
\end{array}
\end{equation*}
and
\begin{equation*}
\begin{array}{ll}
\multicolumn{2}{l}{|D\jump{\varpi}_{r}^{\rm ex} |(\R_N^+) = \ds\int_{\R_+^N} |(\nabla \jump{\varpi}_{r}^{\rm ex})(\xi)| \,d\xi}
\\[2.5ex]
\qquad \qquad & \leq \ds \int_{\R_+^N} |(\tilde{\nabla} \jump{\varpi}_{r}^{\rm ex})(\xi)| \, d\xi + \int_{\R_+^N} |(\partial_N \jump{\varpi}_{r}^{\rm ex})(\xi)| \, d\xi
\\[2.5ex]
& = \ds \int_{\R_+^N}\left|[1 - r^{-1}\xi_N]^+ \tilde{\nabla} \varpi(\tilde{\xi})\right| \, d\xi + \int_{\R_+^N} \left|-r^{-1} \chi_{(0,r)}(\xi_N)\varpi(\tilde{\xi})\right| \, d\xi
\\[2.5ex]
& = \ds\frac{r}{2} \int_{\R^{N-1}}|\tilde{\nabla}\varpi| \, d\tilde{\xi} + |\varpi|_{L^1(\R^{N-1})}
\\[2.5ex]
& = \ds |\varpi|_{L^1(\R^{N-1})} +\tau, \mbox{ for any $ r \in (0, r_{\varpi}^\tau] $.}
\end{array}
\end{equation*}
$\Box$
}

\begin{lemma}\label{LemExt02}
For any $ \hat{v}_\Gamma \in H^1(\Gamma) $ and any $ \ell \in \N $, there exists a function $ \hat{v}_\ell \in  H^1(\Omega) $ such that
\begin{equation}\label{Lem20-2}
\hat{v}_\ell\trace = \hat{v}_\Gamma \mbox{ \ in $ H^{\frac{1}{2}}(\Gamma) $, \ for $ \ell = 1, 2, 3, \dots $,}
\end{equation}
\begin{equation}\label{Lem20-0}
\hat{v}_\ell(x) = 0, \mbox{ \ for a.e. $ x \in \Omega \setminus \Gamma(2^{-\ell}) $ and \ $ \ell = 1, 2, 3, \dots $,}
\end{equation}
and
\begin{equation}\label{Lem20-1}
\ds \left| \hat{v}_\ell \right|_{L^2(\Omega)} \leq 2^{-\ell}
\mbox{ and }
\ds \left|D \hat{v}_\ell \right|(\Omega) \leq |\hat{v}_\Gamma|_{L^1(\Gamma)} +2^{-\ell},
\mbox{ \ for $ \ell = 1,2,3, \dots $.}
\end{equation}
\end{lemma}

\paragraph{Proof of Lemma \ref{LemExt02}.}{
Let $ \sigma > 0 $ be arbitrary, and let $ \rho_*^\sigma $ be the constant as in ($ \Omega $2). Since $ \Gamma \subset \R^{N-1} $ is compact, we can take a large number $ m_\Omega^\sigma \in \N $ and a finite sequence $ \{ x_{\Gamma, 1}^\sigma, \dots, x_{\Gamma, m_\Omega^\sigma}^\sigma \} \subset \Gamma $, such that:
\begin{equation}\label{set01}
\begin{array}{c}
	\ds \overline{\Gamma(r_*/2)} \subset G_*^\sigma := \bigcup_{i = 1}^{m_\Omega^\sigma} G_i^\sigma, \mbox{ \,with the neighborhoods }
\\[3ex]
    G_i^\sigma := G_{x_{\Gamma, i}^\sigma}(\rho_*^\sigma, r_*), ~ i = 1, \dots, m_{\Omega}^\sigma, \mbox{ as in ($\Omega$1)};
\end{array}
\end{equation}
and then, we can take the partition of unity $ \{ \eta_i^\sigma \}_{i = 1}^{m_\Omega^\sigma} \subset C_{\rm c}^\infty(\R^N) $ for the covering $ G_*^\sigma $, such that:
\begin{equation}\label{set02}
0 \leq \eta_i^\sigma \in C_{\rm c}^\infty(G_i^\sigma) \mbox{ \ for $ i = 1, \dots, m_\Omega^\sigma $,\  and \ } \sum_{i = 1}^{m_\Omega^\sigma} \eta_i^\sigma = 1 \mbox{ on $ \overline{\Gamma(r_*/2)} $.}
\end{equation}

Next, let us take any $ \tau > 0 $, and with ($\Omega$1) and Lemma \ref{LemExt01} in mind, let us set:
\begin{equation}\label{set03}
\Xi_i^\sigma := \Xi_{x_{\Gamma, i}^\sigma}, \mbox{ with } \Lambda_i^\sigma := \Lambda_{x_{\Gamma, i}^\sigma} \mbox{ and \ } H_i^\sigma := H_{x_{\Gamma, i}^\sigma}, ~ i = 1, \dots, m_{\Omega}^\sigma,
\end{equation}
\begin{equation}\label{set04}
\begin{array}{c}
\varpi_i^\sigma(\tilde{\xi}) := \left\{ \begin{array}{ll}
\multicolumn{2}{l}{(\eta_i^\sigma \hat{v}_\Gamma) \bigl( (\Xi_i^\sigma)^{-1} \tilde{\xi} \bigr),} 
\\[0.5ex]
& \mbox{ if $ \tilde{\xi} \in \rho_*^\sigma \bB^{N -1} $ and $ i=1,\dots,m_{\Omega}^{\sigma} $,}
\\[1ex]
0, & \mbox{ otherwise,}
\end{array} \right. 
\end{array}
\mbox{for a.e. $ \tilde{\xi} \in \R^{N-1} $,}
\end{equation}
and
\begin{equation}\label{set05}
\hat{r}_{\sigma}^{\tau} := \min \left\{ \begin{array}{l|l}
r_{\varpi_i^\sigma}^{\tau} &  i = 1,\dots,m_\Omega^\sigma
\end{array} \right\}.
\end{equation}

Based on these, we define a class of functions $ \{ \hat{v}_{\sigma}^{\tau} \, | \, \sigma, \tau > 0 \} $, as follows:
\begin{equation}\label{set06}
\begin{array}{c}
\hat{v}_{\sigma}^{\tau}(x) := \left\{ \begin{array}{ll}
\multicolumn{2}{l}{\ds \sum_{i = 1}^{m_\Omega^\sigma} \jump{\varpi_i^\sigma}_{\hat{r}_{\sigma}^{\tau}}^{\rm ex} \bigl( \Xi_i^\sigma x \bigr),} 
\\[0.5ex]
& \mbox{ if $ x \in G_i^\sigma $, for some $ i \in \{1,\dots,m_{\Omega}^{\sigma}\} $,}
\\[1ex]
0, & \mbox{ otherwise,}
\end{array} \right. 
\\[6ex]
\mbox{for a.e. $ x \in \Omega $ and all $ \sigma, \tau > 0 $.}
\end{array}
\end{equation}
Then, as direct consequences of \eqref{set01}--\eqref{set06} and Lemma \ref{LemExt01}, it is inferred that:
\begin{equation}\label{prpty01}
\begin{array}{c}
\hat{v}_\sigma^{\tau} \in H^1(\Omega), ~ \hat{v}_\sigma^{\tau}\trace = \hat{v}_\Gamma \mbox{ in $ H^{\frac{1}{2}}(\Gamma) $,} 
\\[1ex]	
\mbox{ and } \hat{v}_\sigma^{\tau} = 0 \mbox{ a.e. on $ \Omega \setminus \Gamma(\tau) $, for all $ \sigma, \tau > 0 $.}
\end{array}
\end{equation}
Also, in the light of \eqref{Lem10-1}, ($\Omega$2) and Lemma \ref{LemExt01}, we compute that:
\begin{align}
|\hat{v}_\sigma^{\tau}|_{L^2(\Omega)} & = \left[ \rule{-1pt}{18pt} \right. \int_\Omega \Bigl|\sum_{i = 1}^{m_\Omega^\sigma} \jump{\varpi_i^\sigma}_{\hat{r}_\sigma^{\tau}}^{\rm ex}(\Xi_i^\sigma x) \Bigr|^2 \, dx \left. \rule{-1pt}{18pt} \right]^{\frac{1}{2}} \leq \sum_{i = 1}^{m_\Omega^\sigma} \left[ \rule{-1pt}{18pt} \right. \int_{\R_+^{N}}  \bigl| \jump{\varpi_i^\sigma}_{\hat{r}_\sigma^{\tau}}^{\rm ex}(\xi) \bigr|^2 \, d \xi \left. \rule{-1pt}{18pt} \right]^{\frac{1}{2}}\nonumber
\\[1ex]
& \leq m_\Omega^\sigma \tau, \mbox{ for all $ \sigma, \tau > 0 $,}\label{prpty02}
\end{align}
and
\begin{align}
& \hspace{-7ex} \int_\Omega |\nabla_x \hat{v}_\sigma^{\tau}(x)| \, dx \leq \sum_{i =1}^{m_\Omega^\sigma} \int_{G_i^\sigma \cap \Omega} \bigl| \nabla_x \jump{\varpi_i^\sigma}_{\hat{r}_\sigma^{\tau}}^{\rm ex}(\Xi_i^\sigma x) \bigr| \, dx 
\nonumber
\\
& = \sum_{i =1}^{m_\Omega^\sigma} \int_{Y_i^\sigma \cap (\Lambda_i^\sigma \Omega)} \bigl| \nabla_y \jump{\varpi_i^\sigma}_{\hat{r}_\sigma^{\tau}}^{\rm ex}(H_i^\sigma y) \bigr| \, dy
\nonumber
\\
& \leq \sum_{i =1}^{m_\Omega^\sigma} (1 +|\nabla \gamma_{x_{\Gamma}}|_{C(\rho_*^\sigma \overline{\bB^{N -1}})}) \int_{\R_+^N} \bigl| \nabla_{\xi} \jump{\varpi_i^\sigma}_{\hat{r}_{\sigma}^{\tau}}^{\rm ex}(\xi) \bigr| \, d \xi
\nonumber
\\
& \leq (1 +\sigma) \sum_{i = 1}^{m_\Omega^\sigma} \left( \int_{\R^{N -1}} |\varpi_i^\sigma(\tilde{\xi})| \, d \tilde{\xi} +\tau \right)
\nonumber
\\
& \leq (1 +\sigma) \sum_{i = 1}^{m_\Omega^\sigma} \left( \int_{G_i^\sigma \cap \Gamma} \eta_i^\sigma |\hat{v}_\Gamma| \, d \Gamma +\tau \right)
\nonumber
\\
& \leq (1 +\sigma) |\hat{v}_\Gamma|_{L^1(\Gamma)} +m_\Omega^\sigma \tau (1 +\sigma),
\mbox{ for all $ \sigma, \tau > 0 $.}
\label{prpty03}
\end{align}

Now, for any $ \ell \in \N $, let us take two constants $ \sigma_\ell, \tau_\ell \in (0,1] $, such that:
\begin{equation}\label{prpty04}
\left\{ \begin{array}{l}
 ~ (1 +\sigma_\ell) |\hat{v}_\Gamma|_{L^1(\Gamma)} \leq |\hat{v}_\Gamma|_{L^1(\Gamma)} +2^{-\ell-1},
\\[2ex]
 ~ \tau_\ell + m_\Omega^{\sigma_\ell} \tau_\ell (1 + \sigma_\ell) \leq 2^{-\ell-1},
\end{array} \right. \mbox{ for $ \ell = 1, 2, 3, \dots $.}
\end{equation}
Then, on account of \eqref{prpty01}--\eqref{prpty04}, we will conclude that the function $ \hat{v}_\ell := \hat{v}_{\sigma_\ell}^{\tau_\ell} \in H^1(\Omega) $, for each $\ell \in \N $, will fulfill the required condition \eqref{Lem20-2}--\eqref{Lem20-1}. \hfill $\Box$
}

\paragraph{Proof of Key-Lemma A.}{
The proof is a modified version of \cite[Theorem 6]{Moll05}.
Let $ u \in BV(\Omega) \cap L^2(\Omega) $ be arbitrary. Then, by the smoothness  of $ \Gamma $ as in ($ \Omega $1)--($ \Omega$2), we can apply the standard regularization method of BV-functions (cf. \cite[Theorem 10.1.2]{ABM}), and can find a sequence $ \{\hat{\varphi}_\ell\}_{\ell=1}^{\infty} \subset C^{\infty}(\overline{\Omega}) $, such that:
\begin{equation}\label{conv01}
\hat{\varphi}_\ell \to \hat{u} \mbox{ \ in $ L^2(\Omega) $ and strictly in $ BV(\Omega) $, as $ \ell \to \infty $.}
\end{equation}
Besides, from Remark \ref{Rem.trace}, it follows that:
\begin{equation}\label{conv02}
\hat{\varphi}_\ell\trace \to \hat{u}\trace \mbox{ \ in $ L^1(\Gamma) $, as $ \ell \to \infty $.}
\end{equation}

Next, for any $ \ell \in \N $, we apply Lemma \ref{LemExt02} as the case when $ \hat{v}_\Gamma := \hat{u}_\Gamma - \hat{\varphi}_\ell\trace $ in $ H^{\frac{1}{2}}(\Gamma) $, and then, we can take a function $ \hat{\psi}_\ell \in H^1(\Omega) $, such that:
\begin{equation}\label{conv03}
\left\{ ~ \begin{array}{l}
\ds\hat{\psi}_\ell\trace = \hat{u}_\Gamma - \hat{\varphi}_\ell\trace \mbox{ \ in $ H^{\frac{1}{2}}(\Gamma) $,}
\\[1ex]
\ds|\hat{\psi}_\ell|_{L^2(\Omega)} \leq 2^{-\ell} \mbox{ and } \int_{\Omega} |\nabla \hat{\psi}_\ell | \, dx \leq \int_{\Gamma} \left| \hat{u}_\Gamma - \hat{\varphi}_\ell\trace \right| \, d\Gamma + 2^{-\ell}.
\end{array}\right.
\end{equation}

Based on these, let us define:
\begin{equation}\label{conv04}
\hat{u}_\ell := \hat{\varphi}_\ell + \hat{\psi}_\ell \mbox{ \ in $ L^2(\Omega) $, \ for $ \ell =1,2,3,\dots $.}
\end{equation}
Then, in the light of \eqref{conv01}--\eqref{conv03}, it is computed that:
\begin{equation}\label{conv05}
\hat{u}_\ell\trace = \hat{\varphi}_\ell\trace + \hat{\psi}\trace = \hat{\varphi}_\ell\trace + (\hat{u}_\Gamma - \hat{\varphi}_\ell\trace) = \hat{u}_\Gamma \mbox{ \ in $ H^{\frac{1}{2}}(\Gamma) $, for $ \ell=1,2,3,\dots $,}
\end{equation}
\\[-1cm]
\begin{align}
| \hat{u}_\ell - \hat{u} |_{L^2(\Omega)}& =| (\hat{\varphi}_\ell - \hat{u}) + \hat{\psi}_\ell |_{L^2(\Omega)}\nonumber
\\
& \leq | \hat{\varphi}_\ell - \hat{u}|_{L^2(\Omega)} + 2^{-\ell}\to 0 \mbox{ \ as $ \ell \to \infty $,}\label{conv06}
\end{align}
and
\begin{align}
& \hspace{-7ex} \varlimsup_{\ell \to \infty} \int_{\Omega} |\nabla \hat{u}_\ell| \, dx \leq \lim_{\ell\to\infty} \int_{\Omega} |\nabla \hat{\varphi}_\ell| \, dx + \varlimsup_{\ell \to \infty}\int_{\Omega} |\nabla \hat{\psi}_\ell| \, dx
\nonumber
\\
& \leq \int_{\Omega} |D\hat{u}| + \lim_{\ell \to \infty} \left( \int_{\Gamma} |\hat{u}_\Gamma - \hat{\phi}_\ell\trace | \, d\Gamma + 2^{-\ell} \right)
\nonumber
\\
& = \int_{\Omega} |D\hat{u}| + \int_{\Gamma} |\hat{u}\trace -\hat{u}_\Gamma| \, d\Gamma.
\label{conv07} 
\end{align}
Additionally, having in mind Remark \ref{Rem.Du_Dirichlet} (Fact\,4) and \eqref{conv05}--\eqref{conv06}, one can also see that:
\begin{align}
&\hspace{-7ex} \varliminf_{\ell \to \infty} \int_{\Omega} |\nabla \hat{u}_\ell| \, dx = \varliminf_{\ell \to \infty} \left( \int_{\Omega} |\nabla \hat{u}_\ell| \, dx + \int_{\Gamma} |\hat{u}_\ell\trace - \hat{u}_\Gamma | \, d\Gamma \right)
\nonumber
\\
& \geq \int_{\Omega} |D\hat{u}| + \int_{\Gamma} |\hat{u}\trace -\hat{u}_\Gamma| \, d\Gamma.
\label{conv08} 
\end{align}

On account of \eqref{conv05}--\eqref{conv08}, we conclude that the sequence $ \{ \hat{u}_\ell \}_{\ell=1}^{\infty} \subset H^1(\Omega) $, given by \eqref{conv04}, is the required sequence, fulfilling \eqref{Lem100-1}--\eqref{Lem100-2}. \hfill ~ $\Box$
}

\paragraph{Proof of Key-Lemma\,B.}{Let us set:
\begin{equation*}
\sK_0:= \left\{\begin{array}{l|l} W=[w,w_\Gamma] \in \sH & \parbox{5.5cm}{~$ w \leq 0 $, a.e. in $ \Omega $ and $ w_\Gamma \leq 0 $, a.e. on $ \Gamma $} \end{array} \right\}.
\end{equation*}
Then, by using the orthogonal projection $ \pi_{\sK_0} :\sH \to \sK_0 $, we can reformulate the conclusion \eqref{Lem40-0} to the following equivalent form:
\begin{equation}\label{Lem40-1}
\begin{array}{c}
\bigl( U^{*,1}-U^{*,2},(U^1-U^2)-\pi_{\sK_0}(U^1-U^2) \bigr)_\sH \geq 0, 
\\[1ex]
\mbox{ for all $ [U^k,U^{*,k}] \in \partial \Phi_* $ in $ \sH \times \sH $, $ k=1,2 $.}
\end{array}
\end{equation}
Here, according to the general theory of T-monotonicity \cite{KMN80}, the above \eqref{Lem40-1} is equivalent to:
\begin{align}
&\Phi_*(W^1-\pi_{\sK_0}(W^1-W^2)) +\Phi_*( W^2+\pi_{\sK_0}(W^1-W^2))\nonumber
\\
& \leq \Phi_*(W^1) + \Phi_*(W^2) \mbox{, for all $ W^k \in \sW $, $ k=1,2 $.}\nonumber
\end{align}
Additionally, from the definition of $ \sK_0 $, one can easily check that:
\begin{equation*}
\left\{\begin{array}{l}
W^1 -\pi_{\sK_0}(W^1-W^2) = W^1 \vee W^2,
\\[0.5ex]
W^2 +\pi_{\sK_0}(W^1-W^2) = W^1 \wedge W^2,
\end{array}\right.
\mbox{ for all $ W^k \in \sW $, $ k=1,2 $.}
\end{equation*}

Based on these, our goal can be reduced to the verification of:
\begin{equation}\label{Lem40-3}
\begin{array}{c}
\Phi_*(W^1 \vee W^2) + \Phi_*(W^1 \wedge W^2) \leq \Phi_*(W^1) + \Phi_*(W^2),
\\[0.5ex]
\mbox{ for all $ W^k \in D(\Phi_*) $, $ k=1,2 $.}
\end{array}
\end{equation}
Now, to verify \eqref{Lem40-3}, we apply Key-Lemma A, and we can prepare two sequences $ \{ V_\ell^k =  [v_\ell^k,v_{\Gamma,\ell}^k] \}_{\ell =1}^\infty \subset \sV $, $ k=1,2 $, such that:
\begin{equation}\label{Lem40-4}
v_\ell^k\trace = v_{\Gamma,\ell}^k = w_\Gamma^k \mbox{ \ in $ H^{\frac{1}{2}}(\Gamma) $, for every $ \ell \in \N $ and $ k=1,2 $,}
\end{equation}
\begin{equation}\label{Lem40-5}
\begin{array}{c}
\ds v_\ell^k \to w^k \mbox{ in $ L^2(\Omega) $ and } \int_\Omega |\nabla v_\ell^k|\,dx \to \int_\Omega |Dw^k| + \int_\Gamma |w^k\trace - w_\Gamma^k| \,d\Gamma,
\\[2ex]
\mbox{as $ \ell \to \infty $, for every $ k=1,2 $.}
\end{array}
\end{equation}
Subsequently, we compute that:
\begin{align}
& \Phi_*(V_\ell^1 \vee V_\ell^2) + \Phi_*(V_\ell^1 \wedge V_\ell^2)\nonumber
\\
= & \int_\Omega |\nabla v_\ell^1|\,dx + \int_\Omega |\nabla v_\ell^2|\,dx + \frac{\varepsilon^2}{2}\int_\Gamma |\nabla_\Gamma w_\Gamma^1|^2\, d\Gamma + \frac{\varepsilon^2}{2}\int_\Gamma |\nabla_\Gamma w_\Gamma^2|^2\, d\Gamma,\nonumber
\\
& \hspace{3cm}\mbox{ for any $ \ell \in \N $.}\nonumber
\end{align}

Now, taking into account \eqref{Lem40-4}--\eqref{Lem40-5} and the convergences:
\begin{equation*}
V_\ell^1 \vee V_\ell^2 \to W^1 \vee W^2 \mbox{ and } V_\ell^1 \wedge V_\ell^2 \to W^1 \wedge W^2 \mbox{ in $ \sH $ as $ \ell \to \infty $,}
\end{equation*}
the inequality \eqref{Lem40-3} is deduced as follows:
\begin{align*}
&\Phi_*(W^1 \vee W^2) + \Phi_*(W^1 \wedge W^2)
\\
& \leq \lim_{\ell\to \infty} \left( \int_\Omega |\nabla v_\ell^1|\,dx + \int_\Omega |\nabla v_\ell^2|\,dx \right) + \frac{\varepsilon^2}{2}\int_\Gamma |\nabla_\Gamma w_\Gamma^1|^2\, d\Gamma + \frac{\varepsilon^2}{2}\int_\Gamma |\nabla_\Gamma w_\Gamma^2|^2\, d\Gamma
\\
& = \Phi_*(W^1) + \Phi_*(W^2).
\end{align*}
$ \Box $
}

\section{Proofs of the results}

In this section, we prove the results by means of the lemmas and remarks prepared in previous sections.

\paragraph{Proof of Main Theorem A.}
We begin with the verification of the part of lower-bound condition of Mosco-convergence.

Let us take any $ \check{W} = [\check{w}, \check{w}_{\Gamma}] \in \sH $ and any sequence $ \{\check{W}_n = [\check{w}_n,\check{w}_{\Gamma,n}] \}_{n=1}^{\infty} \subset \sV $ such that
\begin{equation*}
\check{W}_n=[\check{w}_n,\check{w}_{\Gamma,n}] \to \check{W} = [\check{w},\check{w}_{\Gamma}] \mbox{ weakly in $ \sH $, as $ n \to \infty $.}
\end{equation*}
Then, for the verification of the inequality of lower-bound condition:
\begin{equation}\label{prMT01}
\varliminf_{n\to\infty} \Phi_n(\check{W}_n) \geq \Phi_*(\check{W}),
\end{equation}
the situation can be restricted to the case that:
\begin{equation}\label{prMT02}
\begin{array}{c}
\ds\lim_{\ell \to \infty} \Phi_{n_\ell}(\check{V}_\ell) = \varliminf_{n\to\infty} \Phi_n(\check{W}_n) < \infty,
\mbox{ \ for some subsequences $ \{ n_\ell \}_{\ell=1}^{\infty} \subset \{ n \} $,}
\\[2ex]
\mbox{ and }\{ \check{V}_\ell =[\check{v}_\ell,\check{v}_{\Gamma,\ell}] \}_{\ell=1}^{\infty} := \{ \check{W}_{n_\ell} =[\check{w}_{n_\ell},\check{w}_{\Gamma,n_\ell}] \}_{\ell=1}^{\infty} \subset \{ \check{W}_n \},
\end{array}
\end{equation}
because the other ones can be said as trivial. Also, from (A1), we can see that:
\begin{align}\label{prMT03}
\Phi_*(\check{V}_{\ell}) & \leq \int_{\Omega} |\nabla \check{v}_\ell| \, dx + \frac{\kappa_{n_\ell}^2}{2} \int_{\Omega} |\nabla \check{v}_\ell|^2 \, dx + \frac{\varepsilon^2}{2} \int_{\Gamma} |\nabla_{\Gamma} \check{v}_{\Gamma,\ell}|^2 \, d\Gamma \nonumber
\\[1ex]
& \leq \Phi_{n_\ell}(\check{V}_\ell) + \mathcal{L}^N(\Omega) \sup_{\omega\in\R^N}\bigl| f_{\delta_{n_\ell}}(\omega) - |\omega|\bigr|, \mbox{ \ for $ \ell=1,2,3,\dots $.}
\end{align}
The conditions \eqref{prMT02}--\eqref{prMT03} imply the boundedness of the sequence $ \{ \check{V}_\ell\}_{\ell=1}^{\infty} $ ($ \subset \sV $) in $ \sW $, and in addition, the assumption (A1) and the lower semi-continuity of $ \Phi_* $ leads to the inequality \eqref{prMT01} of lower-bound condition, via the following calculation:
\begin{equation*}
\begin{array}{ll}
\multicolumn{2}{l}{\ds \varliminf_{n \to \infty} \Phi_n(\check{W}_n) = \ds\lim_{\ell\to\infty} \Phi_{n_\ell}(\check{V}_\ell)}
\\[2ex]
 & \hspace{0.5cm}\geq \ds\varliminf_{\ell\to\infty} \Phi_*(\check{V}_\ell) - \L^N(\Omega) \lim_{\ell\to\infty}\sup_{\omega\in\R^N}\bigl|f_{\delta_{n_\ell}}(\omega) - |\omega| \bigr| \geq  \Phi_*(\check{W}).
\end{array}
\end{equation*}

Next, we show the part of optimality condition. This part can be obtained by applying (A1), Key-Lemma A and the diagonal argument. 

Let us fix any function $ \hat{W} = [\hat{w},\hat{w}_{\Gamma}] \in \sW $. Then, Key-Lemma A enables us to take a sequence $ \{ \hat{V}_\ell =[\hat{v}_\ell,\hat{v}_{\Gamma,\ell}]\}_{\ell=1}^{\infty} \subset \sV $, such that:
\begin{equation}\label{prMT10}
\hat{v}_\ell\trace = \hat{v}_{\Gamma,\ell} = \hat{w}_{\Gamma} \mbox{ \ in $ H^{\frac{1}{2}}(\Gamma) $, \ for $ \ell=1,2,3,\dots $,}
\end{equation}
\begin{equation}\label{prMT11}
\left\{
\begin{array}{l}
\ds \hat{V}_\ell = [\hat{v}_\ell,\hat{v}_{\Gamma,\ell}] \to \hat{W} = [\hat{w},\hat{w}_\Gamma] \mbox{ in $ \sH $,}
\\[2ex]
\ds\int_{\Omega} |\nabla \hat{v}_\ell| \, dx \to \int_{\Omega} |D\hat{w}| + \int_{\Gamma} |\hat{w}\trace - \hat{w}_{\Gamma}| \, d\Gamma,
\end{array}\right. \mbox{ as $ \ell \to \infty $.}
\end{equation}
Here, for any $ \ell \in\N $, let us take a large number $ \hat{n}_\ell \in \N $ such that:
\begin{equation}\label{prMT12}
\frac{\kappa_n^2}{2} \int_{\Omega} |\nabla \hat{v}_\ell|^2 \, dx \leq 2^{-\ell}, \mbox{ \ for any $ n \geq \hat{n}_\ell $.}
\end{equation}
Besides, we define a sequence $ \{ \hat{W}_n=[\hat{w}_n,\hat{w}_{\Gamma,n}] \}_{n=1}^{\infty} \subset \sV $, by letting
\begin{equation}\label{prMT13}
\hat{W}_n = [\hat{w}_n,\hat{w}_{\Gamma,n}] := 
\left\{\begin{array}{ll}
\multicolumn{2}{l}{\hat{V}_\ell = [\hat{v}_\ell,\hat{v}_{\Gamma,\ell}] \mbox{ \ in $ \sV $,}} 
\\[1ex]
& \mbox{ if $ \hat{n}_\ell \leq n < \hat{n}_{\ell+1} $, for some $ \ell \in \N $,}
\\[2ex]
\multicolumn{2}{l}{\hat{V}_1 = [\hat{v}_1,\hat{v}_{\Gamma,1}] \mbox{ \ in $ \sV $,}}
\\[1ex]
& \mbox{ if $ 1 \leq n < \hat{n}_1 $,}
\end{array}\right.
\mbox{ \ for $ n=1,2,3,\dots $.}
\end{equation}
Then, on account of the \eqref{prMT10}--\eqref{prMT13}, it is inferred that
\begin{equation*}
\begin{array}{lll}
\multicolumn{3}{l}{\ds\left|\Phi_n(\hat{W}_n) - \Phi_*(\hat{W}) \right|}
\\[2ex]
& \leq & \ds \left|\int_{\Omega} \left(f_{\delta_n}(\nabla \hat{w}_n) + \frac{\kappa_n^2}{2} |\nabla \hat{w}_n|^2 \right) \, dx - \left(\int_{\Omega} |D\hat{w}| + \int_{\Gamma} |\hat{w}\trace - \hat{w}_{\Gamma}| \,d\Gamma \right)\right|
\\[2ex]
& & \qquad \ds + \frac{\varepsilon^2}{2} \left|\int_{\Gamma} \left(| \nabla_{\Gamma} \hat{w}_{\Gamma,n}|^2 - |\nabla_{\Gamma} \hat{w}_{\Gamma}|^2 \right) \,d\Gamma \right|
\\[2ex]
& \leq & \ds\left|\int_{\Omega} |\nabla \hat{w}_n| \, dx - \left(\int_{\Omega} |D\hat{w}| + \int_{\Gamma} |\hat{w}\trace - \hat{w}_{\Gamma} | \, d\Gamma \right)\right|
\\[3ex]
& &\qquad \ds +\mathcal{L}^N(\Omega) \sup_{\omega\in\R^N}\bigl| f_{\delta_n}(\omega) - |\omega| \bigr| + 2^{-\ell},
\\[3ex]
\multicolumn{3}{c}{\mbox{for any $ \hat{n}_\ell \leq n < \hat{n}_{\ell+1} $, $ \ell =1,2,3,\dots $,}}
\end{array}
\end{equation*}
and it implies the convergence $ \lim_{n\to\infty} \Phi_n(\hat{W}_n) = \Phi_*(\hat{W}) $, required in optimality condition.

Thus, we conclude Main Theorem A. \hfill $ \Box $
\begin{remark}\label{remMainA}
\begin{em}
Let us simply denote by $ \Phi_0 := \Phi_* |_\sV $ the restriction of $ \Phi_* $ onto $ \sV $, more precisely:
\begin{center}
$ \ds V = [v, v_\Gamma] \in \sV \mapsto \Phi_0(V) = \Phi_0(v, v_\Gamma) := \int_{\Omega} |\nabla v| \, dx +\frac{\varepsilon^2}{2} \int_\Gamma |\nabla_{\Gamma} v_{\Gamma}|^2 \, d\Gamma $.
\end{center}
Then, as a consequence of Main Theorem A, one can observe that $ \Phi_* $ coincides with the lower semi-continuous envelope $ \overline{\Phi_0} $ of the restriction $ \Phi_0 $, i.e.:
\begin{equation}\label{l.s.c.env}
\begin{array}{c}
\Phi_*(W) = \overline{\Phi_0}(W) := \inf \left\{ \begin{array}{l|l}
\ds \varliminf_{n \to \infty} \Phi_0(V_n) & 
 ~ \parbox{4.75cm}{
$ \{ V_n \}_{n = 1}^\infty \subset \sV $ \ and \\[0.25ex] $ V_n \to W $ in $ \sH $ as $ n \to \infty $
}
\end{array} \right\},
\\
\ \\[-2ex]
\mbox{for any $ W \in \mathscr{H} $.}
\end{array}
\end{equation}
In fact, from \eqref{l.s.c.env} of $ \overline{\Phi_0} $, we see that the lower semi-continuous envelope $ \overline{\Phi_0} $ is a maximal l.s.c. function supporting $ \Phi_0 $ on $ \sV $. So, we immediately have:
\begin{equation}\label{cor01}
\Phi_* \leq \overline{\Phi_0} \mbox{ \ on $ \sH $, and } D(\overline{\Phi_0}) \subset D(\Phi_*) = \sW. 
\end{equation}
Meanwhile, for any $ \hat{W} = [\hat{w},\hat{w}_{\Gamma}] \in D(\overline{\Phi_0}) $, taking the sequence $ \{ \hat{V}_\ell =[\hat{v}_\ell, \hat{v}_{\Gamma,\ell}] \}_{\ell=1}^{\infty} \subset \sV $, as in \eqref{prMT10}--\eqref{prMT11}, enables us to deduce that:
\begin{equation}\label{cor02}
\overline{\Phi_0}(\hat{W}) \leq \lim_{\ell \to \infty} \Phi_0(\hat{V}_\ell) = \Phi_*(\hat{W}).
\end{equation} 

\eqref{cor01} and \eqref{cor02} imply the coincidence $ \Phi_*=\overline{\Phi_0} $ on $ \sH $.
\end{em}
\end{remark}
\paragraph{Proof of Corollary \ref{Cor.03}}
 This corollary will be obtained as straightforward consequences of Main Theorem A and the general theories of abstract evolution equations and their variational convergences, e.g. \cite{Attouch, Barbu, Brezis, Kenmochi81}, and so on.\hfill $ \Box $

\paragraph{Proof of Main Theorem\,B.}{
By the assumption, we find two functions \\ $ U^{*,k} \in L^2(0,T;\sH) $, $ k=1,2 $, such that:
\begin{equation}\label{prMT200}
\begin{array}{c}
U^{*,k}(t) \in \partial \Phi_*(U^k(t)) \mbox{ and }
(U^k)'(t) + U^{*,k}(t) = \Theta^k(t) \mbox{ in $ \sH $,}
\\[1ex]
\mbox{for a.e. $ t \in (0,T) $, $ k=1,2 $.}
\end{array}
\end{equation}
Here, taking the difference between the equations in \eqref{prMT200} and multiplying the both sides by $ [U^1 -U^2]^+(t) $, one can see that:
\begin{align}
&\frac{1}{2}\frac{d}{dt}\bigl| [U^1-U^2]^+(t)\bigr|_\sH^2 + \bigl((U^{*,1}-U^{*,2})(t), [U^1 - U^2]^+(t)\bigr)_\sH\nonumber
\\
&\hspace{2cm}=\bigl((\Theta^1-\Theta^2)(t), [U^1-U^2]^+(t) \bigr)_\sH, \mbox{ \ a.e. $ t \in (0,T) $. }\label{prMT201}
\end{align}
Also, from Key-Lemma\,B, it immediately follows that:
\begin{equation}\label{prMT202}
\bigl((U^{*,1}-U^{*,2})(t), [U^1 - U^2]^+(t)\bigr)_\sH \geq 0.
\end{equation}
Thus, Main Theorem B will be concluded by using the standard method, i.e. by applying \eqref{prMT202}, Young's inequality and Gronwall's lemma to \eqref{prMT201}. \hfill $ \Box $
}

\begin{remark}\label{Rem.hosoku}
\begin{em}
In the proofs of Main Theorems A and B, the essentials will be in the fixed-situations of boundary data for approximating functions, as in \eqref{Lem100-1}, \eqref{Lem40-4} and \eqref{prMT10}. 
Then, the auxiliary Lemmas \ref{LemExt01}--\ref{LemExt02} are to support the presence of such approximations, and proofs of these can be said as some simplified version of the regularization method developed by Gagliardo \cite{Gagliardo}. 
But, the original method by \cite{Gagliardo} would be available just for the regularizations of BV-functions by $ W^{1,1} $-functions, and it would not support the regularizations by other kinds of functions, so immediately. 
Hence, for the $ H^1 $-regularizations required in this study, the simplified construction \eqref{Lem10-0}--\eqref{Lem10-1} would be essential, and then, the $ H^1 $-regularity of the boundary data would be needed to be the assumptions, as in Key-Lemma A and Lemmas \ref{LemExt01}--\ref{LemExt02}. 
\end{em}
\end{remark}
%


\section{Future prospective}

One of the possible prospectives is to apply our theory to the phase-field system of grain boundary motion, known as ``Kobayashi--Warren--Carter model'', cf. \cite{KWC1,KWC2}. Indeed, the Kobayashi--Warren--Carter model is derived as a gradient system of a governing energy, including a generalized (unknown-dependent) total variation. In this light, the objective of this issue will be in the enhancement of the mathematical method for grain boundary phenomena, if we can combine our results and the line of relevant works to the Kobayashi--Warren--Carter model, e.g. \cite{IKY08, IKY09, KWC1,KWC2, MS14, SWY13, SWY14}.

\bigskip

\paragraph{Acknowledgments}{
On a final note, we appreciate very much to the anonymous referee for taking great efforts to review our manuscript, and for giving us a lot of valuable comments and remarks. 
}


\begin{thebibliography}{99}

\bibitem{AFP}
     Ambrosio, L.; Fusco, N.; Pallara, D.:
     \newblock \emph{Functions of Bounded Variation and Free Discontinuity Problems.}
     \newblock Oxford University Press,
     \newblock New York (2006).
%
\bibitem{ABCM} 
    Andreu, F.; Ballester, C.; Caselles, V.; Maz\'{o}n, J. M.:
    \newblock The Dirichlet problem for the total variation flow.
    \newblock J. Funct. Anal. \textbf{180} (2001), no. 2, 347--403. 
%
\bibitem{Anzellotti} 
    Anzellotti, G.:
    \newblock The Euler equation for functionals with linear growth.
    \newblock Trans. Amer. Math. Soc. \textbf{290} (1985), 483--501.
%
\bibitem{Attouch}
     Attouch, H.:\emph{Variational Convergence for Functions and Operators.}
     \newblock Applicable Mathematics Series,
     \newblock Pitman, Massachusetts (1984).

\bibitem{ABM}
	\newblock Attouch, H.; Buttazzo, G.; Michaille, G.:  
	\newblock \emph{Variational Analysis in Sobolev and BV Spaces.} 
	\newblock Applications to PDEs and Optimization. 
	\newblock MPS-SIAM Series on Optimization, 6. SIAM and MPS, (2006).

\bibitem{Barbu}
	Barbu, V.:
\newblock \emph{Nonlinear Differential Equations of Monotone Type in Banach Spaces.}
	\newblock Springer Monographs in Mathematics. Springer
	\newblock Springer, New York (2010).%
	
\bibitem{Brezis}
	 Br\'{e}zis, H.: 
	\newblock \emph{Operateurs Maximaux Monotones et Semi-groupes de Contractions dans les Espaces de Hilbert.}
	\newblock North-Holland Mathematics Studies, \textbf{5}, 
	\newblock Notas de Matem\'{a}tica (50), 
	\newblock North-Holland Publishing and American Elsevier Publishing (1973).
	
\bibitem{CC13}
     Calatroni, L.; Colli, P.:
    \newblock Global solution to the Allen-Cahn equation with singular potentials and dynamic boundary conditions.
    \newblock Nonlinear Anal. \textbf{79} (2013), 12--27.

\bibitem{CGNS1X}
    \newblock Colli, P.; Gilardi, G.; Nakayashiki, R.; Shirakawa, K.:
    \newblock A class of quasi-linear Allen--Cahn type equations with dynamic boundary conditions.
    \newblock Nonlinear Anal. \textbf{158} (2017), 32--59.
%
\bibitem{CF1}
	 Colli, P.; Fukao, T.: 
	\newblock The Allen--Cahn equation with dynamic boundary conditions and mass constraints.
    \newblock Math. Methods Appl. Sci. {\bf 38} (2015), 3950--3967.

\bibitem{CS}
    Colli, P.; Sprekels, J.: 
    \newblock Optimal control of an Allen--Cahn equation with singular potentials and dynamic boundary condition.
    \newblock SIAM J. Control Optim. {\bf 53} (2015), 213--234.

\bibitem{EG}
	\newblock Evans, L. C.; Gariepy, R. F.:
	\newblock \emph{Measure Theory and Fine Properties of Functions. Revised edition.}
	\newblock Textbooks in Mathematics, CRC Press, Inc., Boca Raton (2015).

\bibitem{Gagliardo}
   \newblock Gagliardo, E.:
   \newblock Caratterizzazioni delle tracce sulla frontiera relative ad alcune classi di funzioni in n variabili. 
   \newblock Rend. Sem. Mat. Univ. Padova (Italian) \textbf{27} (1957), 284--305.
    
\bibitem{GGM08}
    Gal, C. G.; Grasselli, M.; Miranville, A.:
    \newblock Nonisothermal Allen-Cahn equations with coupled dynamic boundary conditions.
    \newblock Nonlinear phenomena with energy dissipation, 117--139, GAKUTO Internat. Ser. Math. Sci. Appl., 29, Gakkt\={o}sho, Tokyo, 2008.
%
\bibitem{G}
	\newblock Giusti, E.:
	\newblock Minimal Surfaces and Functions of Bounded Variation.
	\newblock Monographs in Mathematics \textbf{80}, Birkh\"{a}user (1984).
%
\bibitem{IKY08}
	\newblock Ito, A.; Kenmochi, N.; Yamazaki, N.: 
	\newblock A phase-field model of grain boundary motion. 
	\newblock Appl. Math., \textbf{53} (2008), no. 5, 433--454.

\bibitem{IKY09}
	\newblock Ito, A.; Kenmochi, N.; Yamazaki, N.: 
	\newblock Weak solutions of grain boundary motion model with singularity. 
	\newblock Rend. Mat. Appl. (7), \textbf{29} (2009), no. 1, 51--63.

\bibitem{KWC1}
	\newblock Kobayashi, R.; Warren, J. A.; Carter, W. C.: 
	\newblock A continuum model of grain boundaries. 
	\newblock Phys. D, \textbf{140} (2000), no. 1-2, 141--150.

\bibitem{KWC2}
	\newblock Kobayashi, R.; Warren, J. A.; Carter, W. C.: 
	\newblock Grain boundary model and singular diffusivity. 
	\newblock In: \emph{Free Boundary Problems: Theory and Applications}, pp. 283--294, 
	\newblock GAKUTO Internat. Ser. Math. Sci. Appl., \textbf{14}, Gakk\={o}tosho, Tokyo, (2000).
%
\bibitem{Kenmochi81}%
	 Kenmochi, N.: 
    \newblock Solvability of nonlinear evolution equations with time-dependent constraints and applications. 
    \newblock Bull. Fac. Education, Chiba Univ., \textbf{30} (1981), 1--87. \hfill\break
    \newblock \url{http://ci.nii.ac.jp/naid/110004715232}

\bibitem{KMN80}%
	\newblock Kenmochi, N.; Mizuta, Y.; Nagai, T.:
	\newblock Projections onto convex sets, convex functions and their subdifferentials.
    \newblock Bull. Fac. Education, Chiba Univ., \textbf{29} (1980), 11--22. \hfill\break
    \newblock \url{http://ci.nii.ac.jp/naid/110004715212}

\bibitem{Moll05}
    Moll, J. S.:
    \newblock The anisotropic total variation flow. 
    \newblock Math. Ann. \textbf{332} (2005), no. 1, 177--218. 
%
\bibitem{MS14} 
	\newblock Moll, S.; Shirakawa, K.: 
	\newblock Existence of solutions to the Kobayashi-Warren-Carter \linebreak system. 
	\newblock Calc. Var. Partial Differential Equations, \textbf{51} (2014), 621--656.  
	\newblock {DOI:10.1007/ s00526-013-0689-2}

\bibitem{Mosco}
	Mosco, U.:
	\newblock Convergence of convex sets and of solutions of variational inequalities.
	\newblock Advances in Math. \textbf{3}, 510--585 (1969).
%
\bibitem{SV97}
	 Savar\'{e}, G.; Visintin, A.:
	\newblock Variational convergence of nonlinear diffusion equations: applications to concentrated capacity problems with change of phase. 
	\newblock Atti Accad. Naz. Lincei Cl. Sci. Fis. Mat. Natur. Rend. Lincei (9) Mat. Appl. \textbf{8}, (1997), no. 1, 49--89. 
%
\bibitem{SWY13} 
    \newblock Shirakawa, K.; Watanabe, H.; Yamazaki, N.: 
    \newblock Solvability of one-dimensional phase field systems associated with grain boundary motion.
    \newblock Math. Ann., \textbf{356} (2013), 301--330. {DOI:10.1007/s00208-012-0849-2}

\bibitem{SWY14}
	Shirakawa, K., Watanabe, H., Yamazaki, N.:
	\newblock Phase-field systems for grain boundary motions under isothermal solidifications.
	\newblock Adv. Math. Sci. Appl., \textbf{24} (2014), 353--400.
%
\bibitem{Temam} 
   Temam, R.:
   \newblock On the continuity of the trace of vector functions with bounded deformation.
   \newblock Appl. Anal. \textbf{11} (1981), 291--302.
\end{thebibliography}
\end{document}